\documentclass[journal]{IEEEtran}

\pdfoutput=1 
\usepackage{amssymb,epsfig,color,cite,amsmath,amsfonts,mathrsfs,algorithm,algorithmicx}
\usepackage{epstopdf}
\usepackage{datetime,fancyhdr}
\usepackage{algpseudocode}
\usepackage{arydshln}
\usepackage{multirow}
\usepackage{float}

\usepackage{times} 
\usepackage{epsfig}
\usepackage{graphicx}
\usepackage{latexsym}
\usepackage{subfigure}

\newtheorem{lemma}{Lemma}[section]
\newtheorem{theorem}{Theorem}[section]

\newtheorem{remark}{Remark}[section]
\newtheorem{definition}{Definition}[section]
\newtheorem{example}{Example}[section]

\newtheorem{assumption}{Assumption}[section]

\newcommand{\upcite}[1]{\textsuperscript{\textsuperscript{\cite{#1}}}}

\ifCLASSINFOpdf
\else
\fi
\begin{document}
%
\title{Distributed synchronous and asynchronous algorithms for semi-definite programming with diagonal constraints}
%
%
%

\author{Xia Jiang,
        Xianlin~Zeng,~\IEEEmembership{Member,~IEEE,}
        Jian~Sun,~\IEEEmembership{Member,~IEEE,}
        and~Jie~Chen,~\IEEEmembership{Fellow,~IEEE}
\thanks{This work was supported in part by the National Natural Science Foundation of China (Nos. 61720106011, U1613225, 62073035, 61925303), Program for Changjiang Scholars and Innovative Research Team in University (IRT1208), the Youth ChangJiang Scholars Program, the National Key Research and Development Program of China under Grant 2018YFB1700100 and Beijing Institute of Technology Research Fund Program for Young Scholars. \emph{(Corresponding author: Jian Sun.)}}
\thanks{X. Jiang (jiangxia@bit.edu.cn) and J. Sun (sunjian@bit.edu.cn) are with Key Laboratory of Intelligent Control and Decision of Complex Systems, School of Automation, Beijing Institute of Technology, Beijing, 100081, China, and also with the Beijing Institute of Technology Chongqing Innovation Center, Chongqing  401120, China}
\thanks{X. Zeng (xianlin.zeng@bit.edu.cn) is with Key Laboratory of Intelligent Control and Decision of Complex Systems, School of Automation, Beijing Institute of Technology, Beijing, 100081, China}
\thanks{J. Chen (chenjie@bit.edu.cn) is with Beijing Advanced Innovation Center for Intelligent Robots and Systems (Beijing Institute of Technology), Key Laboratory of Biomimetic Robots and Systems (Beijing Institute of Technology), Ministry of Education, Beijing, 100081, China, and also with the School of Electronic and Information Engineering, Tongji University, Shanghai, 200082, China}}

\markboth{IEEE Transactions on Automatic Control,~\today~\currenttime}
{\MakeLowercase{\textit{et al.}}: Distribute method}


\maketitle

\begin{abstract}
This paper develops distributed synchronous and asynchronous algorithms for the large-scale semi-definite programming with diagonal constraints, which has wide applications in combination optimization, image processing and community detection. The information of the semi-definite programming is allocated to multiple interconnected agents such that  each agent  aims to find a solution  by communicating to its neighbors. Based on low-rank  property of solutions and  the Burer-Monteiro factorization,  we  transform the original  problem into a distributed  optimization problem  over unit spheres to reduce variable dimensions and ensure positive semi-definiteness without involving semi-definite projections, which are computationally expensive. 
For  the distributed optimization problem, we propose distributed synchronous and asynchronous algorithms, both of which  reduce computational burden and storage space compared with existing centralized algorithms. 
Specifically, the  distributed synchronous algorithm almost surely  escapes strict saddle points and converges to the set of optimal solutions  to the optimization problem. In addition, the proposed distributed asynchronous algorithm allows communication delays and converges to the set of  critical points to the optimization problem under mild conditions. By applying proposed algorithms to image segmentation applications, we illustrate the efficiency and convergence performance of the two proposed algorithms.
\end{abstract}

\begin{IEEEkeywords}
 Semi-definite programming with diagonal constraints, synchronous and asynchronous algorithms, low-rank matrices, distributed optimization.
\end{IEEEkeywords}

%
\IEEEpeerreviewmaketitle

\section{Introduction}
%
%
%
%
\IEEEPARstart{S}{emi-definite} programming (SDP) is an active subfield of convex optimization and has attracted considerable attention due to its widely applications in diverse fields such as control theory\cite{smart_admm,semi_observer,2003Semidefinite}, combinatorial optimization \cite{center_sdp,combina}, operations research\cite{max_like,eco_dispatch}, and machine learning\cite{matrix_learning,kernel_matrixlearning,distance_matrix}. Formally, it aims to maximize or minimize a linear objective function subject to a constraint that is an affine combination of positive semi-definite matrices. One important class of SDP problem is the  SDP with diagonal constraints, which is a relaxation of the ``maximum cut" problem \cite{1995Improved} and also appears in phase retrieval \cite{Waldspurger2012Phase} and $\mathbb Z/2\mathbb Z$ synchronization \cite{2014Exact}. 

Various algorithms have been developed to solve SDP with diagonal constraints but tend to be computational demanding as variable dimensions scale.
On one hand, the arithmetic cost scales badly as the dimension of matrices increases, especially for high-order algorithms. For example, each iteration  costs $O(n^3)$ arithmetic operations with an interior-point solver, which solves SDP in polynomial time \cite{2007Implementation}, and  the computation may run out of memory and time if $n$ is greater than several thousands \cite{BM_smooth}. On the other hand, the storage cost of each iteration may scale beyond the   memory of single computer if the number of unknowns reaches tens of thousands \cite{2015Phase}. Hence, the design of efficient algorithms for large-scale SDP with diagonal constraints is still a challenging problem.

\par  
To reduce the computational burden of large-scale SDP, there are many explorations of efficient centralized works in recent years\upcite{BM_smooth,2017Solving,wang2017mixing,Zhang2012}. 
One key idea is using the Burer-Monteiro factorization that expoits the low-rank property of matrix solutions by replacing the original large scale positive semi-definite matrix  as the product of two ``tall" matrices with lower dimensions to reduce storage cost and avoid computing the semi-definite projection. 
Using this idea, in \cite{BM_SPLR,BM_smooth,lowrank_sdp,2017Solving}, authors transform general SDP into non-convex optimization problems by making use of the low-rank property of solutions, and propose augmented Lagrangian algorithms and Riemannian manifold methods. In addition, the challenging positive semi-definite constraints are eliminated with the cost of introducing non-convexity to the optimization problems. Surprisingly, this change to  non-convex problems does not cause many difficulties because local solutions tend to recover the optimal solution in practice. Despite these advances, the existing augmented Lagrangian algorithms and Riemannian manifold methods do not guarantee converging to global optima, and suffer from slow convergence and difficulties in selecting step sizes. For SDP with diagonal constraints, some recent works\cite{wang2017mixing,erdogdu2018convergence} developed block-coordinate algorithms with rigorous convergence analysis, which have free parameters and better optimization performance than prior works\upcite{BM_SPLR,BM_smooth,lowrank_sdp,2017Solving}. All these centralized algorithms own fine practical evidences for the transformed non-convex optimizations. However, as the matrix dimension grows too large, the lower dimensional matrices in these algorithms may still take too much storage space for a single computer, such as these in some image processing problems. In addition, some information and data of practical problems may be generated and stored at different  locations and cannot be communicated  due to privacy considerations. Hence,  these centralized algorithms can not be applied directly to large scale problems with distributed information and distributed algorithms are in need for large-scale SDP.

\par Distributed optimization algorithms offer a promising approach to address large scale matrix problems by using the problem setup that the information is allocated over different agents\upcite{Sylvester_matrix,linear_eq,Deng2019NetworkFT,2016Implementing,matrix_equation,Lyapunov_matrix,Hong2016}. In distributed setting, agents have access to local information and communicate with their neighbors to seek for a global optimal solution\upcite{consensus_shi,liang_opti,Wang2020}. For large-scale SDP, many works in \cite{smart_admm,ADMM_CDC,fast_dis,semi_power,asyn_admm}  exploited the sparse structure of SDP and introduced additional consensus constrains to the transformed distributed problems. These works proposed distributed algorithms based on alternating direction method of multipliers (ADMM) with iterative message-passing.  Whereas, in ADMM, agents need to solve sub-semidefinite problems at each iteration and have considerable computational burden. Focusing on SDP with tree structures, \cite{dis_pd} proposed a distributed primal-dual interior-point algorithm for constrained semi-definite programming without introducing consensus constraints. The algorithm in \cite{dis_pd} is a second-order algorithm that conducts a recursion over the tree structure to compute the exact search directions and factorizes a relatively small matrix during each iteration. Recently, \cite{ricatti_zeng} proposed a distributed optimization design for solving continuous-time algebraic Riccati inequalities, which have applications in distributed control of multi-agent systems. This design is a first-order algorithm and has well intuitive interpretations, but it needs computing semi-definite projections, which are expensive for large-scale matrices.  



\par In this paper, we  develop distributed first-order algorithms for large-scale semi-definite programs with diagonal constraints by taking advantage of low-rank property of solutions and the inherent sparsity of problems.
 The contributions of this paper are summarized as follows.
\begin{itemize}
\item This paper proposes a study on the distributed algorithms for SDP with diagonal constraints and distributed  coefficient matrices information. This study extends the works in \cite{BM_SPLR,BM_smooth,wang2017mixing,erdogdu2018convergence} to distributed setups, which have wide applications in power flow problems\upcite{ADMM_CDC,smart_admm} and distributed state estimation/control\upcite{semi_power,ricatti_zeng}. In addition, the SDP problem in this paper does not require tree structures as in \cite{dis_pd}.
\item This paper  designs distributed synchronous and asynchronous algorithms for SDP with diagonal constraints by solving an equivalent nonconvex optimization problem, which is obtained using the Burer-Monteiro factorization. 
In particular, the distributed algorithms reduce the computational burden and storage  cost on single agent compared with the existing centralized algorithms\cite{wang2017mixing,erdogdu2018convergence} for SDP with diagonal constraints and show a superior numerical performance in simulation experiments. With the Burer-Monteiro factorization, the proposed algorithms avoid the computational burden of projection to semi-definite cone\upcite{ADMM_CDC,fast_dis,semi_power,asyn_admm}. Compared with the distributed second-order interior-point algorithm in \cite{dis_pd}, the proposed first-order algorithms have lower complexity and the distributed asynchronous algorithm performs well without a global synchronous clock. 
\item This paper analyzes the convergence of our proposed distributed algorithms. For the distributed synchronous algorithm, we show that the variables converge to the set of global optimal solutions almost surely under random initializations, despite of
the non-convexity of feasible sets. For the distributed asynchronous algorithm, we show that the variables converge to the set of critical points of the nonconvex problem under mild conditions.
\end{itemize}
 \par The remainder of the paper is organized as follows. Mathematical notations are given in section \ref{preliminaries_sec}. The semi-definite programming description and distributed algorithms are proposed in section \ref{solver_design}. The convergence properties of the proposed algorithms are analyzed theoretically in section \ref{proof_sec}. The efficiency of distributed algorithms is verified by simulations in section \ref{simulation} and the conclusion is made in section \ref{conclusion}.

\section{Mathematical Notations} \label{preliminaries_sec}

\par We denote $\mathbb{R}$ as the set of real numbers, $\mathbb{R}^n$ as the set of $n$-dimensional real column vectors, $\mathbb{R}^{n\times m}$ as the set of $n$-by-$m$ real matrices, $\mathbb{N}$ as the set of natural numbers, $\mathbb{S}^n$ as the set of $n$ by $n$ symmetric matrices, $\mathbb{S}_+^n$ as the set positive semi-definite matrices, $\emptyset$  as the empty set, respectively. All vectors in the paper are column vectors, unless otherwise noted. The notation $0_n$ denotes an $n \times 1$ vector with all elements of $0$. For a real vector $v$, $\left\|v\right\|$ is the Euclidean norm and $\left\|v\right\|_1$ is 1-norm defined by the sum of  absolute values of elements. We denote $A'$ as the transpose of matrix $A$, $\lambda_{min}(A)$ as the minimum eigenvalue of the matrix $A$. For a symmetric matrix $A$, $A \succeq 0$ denotes that $A$ is positive semi-definite and $A_{(i,j)}$ is the $(i,j)$th element of matrix $A$. For real matrices $A$ and $B$ with same dimensions, $\left<A, B\right>$ denotes the Frobenius inner product of two real matrices such that $\left< A,B\right>=tr(A'B)=\sum_{i,j}A_{(i,j)}B_{(i,j)}$. In addition, $A\circ B$ denotes Hadamard product of two matrices, whose elements are defined by $[A\circ B]_{(i,j)}=A_{(i,j)}B_{(i,j)}$. For a twice-continuously differentiable function $f(x)$, its gradient and Hessian matrix are denoted as $\nabla f(x)$ and $\nabla^2 f(x)$.
\par For a set $\mathcal{S}$, $\left|\mathcal{S}\right|$ denotes the number of elements in the set $\mathcal{S}$. For sets $\mathcal{S}_1$ and $\mathcal{S}_2$, $\mathcal{S}_1\subset \mathcal{S}_2$ means that $\mathcal{S}_1$ is a subset of $\mathcal{S}_2$, $\mathcal{S}_1\cup \mathcal{S}_2$ is the union of $\mathcal{S}_1$ and $\mathcal{S}_2$, $\mathcal{S}_1 \cap \mathcal{S}_2$ is the intersection of $\mathcal{S}_1$ and $\mathcal{S}_2$, and $\mathcal{S}_1\backslash \mathcal{S}_2=\mathcal{S}_1-(\mathcal{S}_1\cap \mathcal{S}_2)$. For a real number $a$, $\lceil a \rceil$ is the smallest integer greater than $a$. For a non-zero vector $x\in \mathbb{R}^n$, the notation ${\rm normal}(x)$ is $\frac{x}{\left\|x\right\|}$ and $A={\rm diag}(x)\in \mathbb{R}^{n\times n}$ denotes a matrix with diagonal element $A_{(i,i)}=x_i$.
\par Let $\mathcal{X}$ be a smooth manifold. Let $f:\mathcal{X} \to \mathbb{R}$ be a real-valued twice-continuously differentiable function. A point $x^*$ is a critical point of $f$ if $\nabla f(x^*)=0_d$. If, in addition, $\lambda_{min}(\nabla^2 f(x^*))<0$, $x^*$ is a strict saddle point of $f$. 
\section{Problem Description and Distributed Algorithm Design}\label{solver_design}

In this section, we present the problem of solving semi-definite programming with diagonal constraints in a  distributed way. Then, we reformulate the problem into a distributed non-convex optimization using the low-rank property of solutions, and propose distributed synchronous and asynchronous discrete-time algorithms.
\subsection{Problem description and transformation}\label{intro_pro_set}

Let $\mathcal{G}=(\mathcal{V},\mathcal{E})$ be an arbitrary simple, undirected and connected graph with the node set $\mathcal{V}=\{1,\cdots,m\}$ and the edge set $\mathcal{E}=\mathcal{V}\times \mathcal{V}$. Let $X\in\mathbb S_+^{n}$ be the variable. For each  $\tilde i\in \mathcal{V}$, define $J_{\tilde{i}}\subset \{1,\cdots,n\}$ as an ordered set, where $n>m$.
Throughout this paper, we use $\tilde{(\cdot)}$ to denote the index of agent in $\mathcal G$ and $(\cdot)$ to denote an element of $J_{\tilde{(\cdot)}}$. For example,  $\tilde i\in \mathcal{V}$ denotes agent $\tilde i$ of $\mathcal G$ and $j\in J_{\tilde{i}}$ is element $j$ in $J_{\tilde{i}}$. If $J_{\tilde{i}}\cap J_{\tilde{j}}\neq \emptyset$, then  $(\tilde{i},\tilde{j})\in \mathcal{E}$  such that  agents $\tilde{i}$ and $\tilde{j}$ can communicate with each other.


The distributed semi-definite programming with diagonal constraints  is
\begin{subequations}\label{sdp_pro}
	\begin{align}
	\min_{X\in\mathbb{S}_+^n} \,\, & \sum_{\tilde{i}=1}^m   \langle\overline{M}^{\tilde{i}},  X\rangle, \label{sdp_func}\\
	{\rm s. \, t.}\, \, & X_{(j,j)}=1, \quad j \in \{1,\cdots,n\}, \label{sdp_constr}
	\end{align}
\end{subequations}
where $\overline{M}^{\tilde{i}}\in \mathbb{S}^n$ is a coefficient matrix such that $\overline{M}^{\tilde{i}}_{(j,k)} = 0$ if $(j,k)\notin J_{\tilde{i}}\times J_{\tilde{i}}$, {and $\sum_{\tilde{i}=1}^m \overline{M}^{\tilde{i}}=M$.} 
Define $E_{J_{\tilde{i}}\times J_{\tilde{i}}}\in \mathbb{R}^{n\times n}$ as the $0-1$ matrix with $E_{(l,k)}=1$ for $(l,k) \in J_{\tilde{i}}\times J_{\tilde{i}}$ and $E_{(l,k)}=0$ otherwise. It is clear that $ \langle\overline{M}^{\tilde{i}},  X\rangle =  \langle\overline{M}^{\tilde{i}},  E_{J_{\tilde{i}}\times J_{\tilde{i}}} \circ X\rangle$ for all $i\in\{1,\ldots,n\}$. Without affecting solutions, the local  variable $\overline{X}^{\tilde{i}}\in\mathbb{S}_{+}^{n}$ to be determined by agent $\tilde{i}\in\mathcal V$ is defined as $\overline{X}^{\tilde{i}}=E_{J_{\tilde{i}}\times J_{\tilde{i}}} \circ X$. 
The objective of this paper is to design a distributed algorithm for solving \eqref{sdp_pro} such that each agent $\tilde{i}\in\mathcal V$ only knows local information $\overline{M}^{\tilde{i}}\in \mathbb{S}^n$.

\begin{example}
	Figure \ref{graph_cut} illustrates the  relationship of global optimization variable $X$ and  local variables $\overline{X}^{\tilde{i}}$ ($i\in \{1,\cdots,4\}$), where elements with same indices in different colored matrices are the same.
	\begin{figure*}
		\centering
		\includegraphics[scale=0.4]{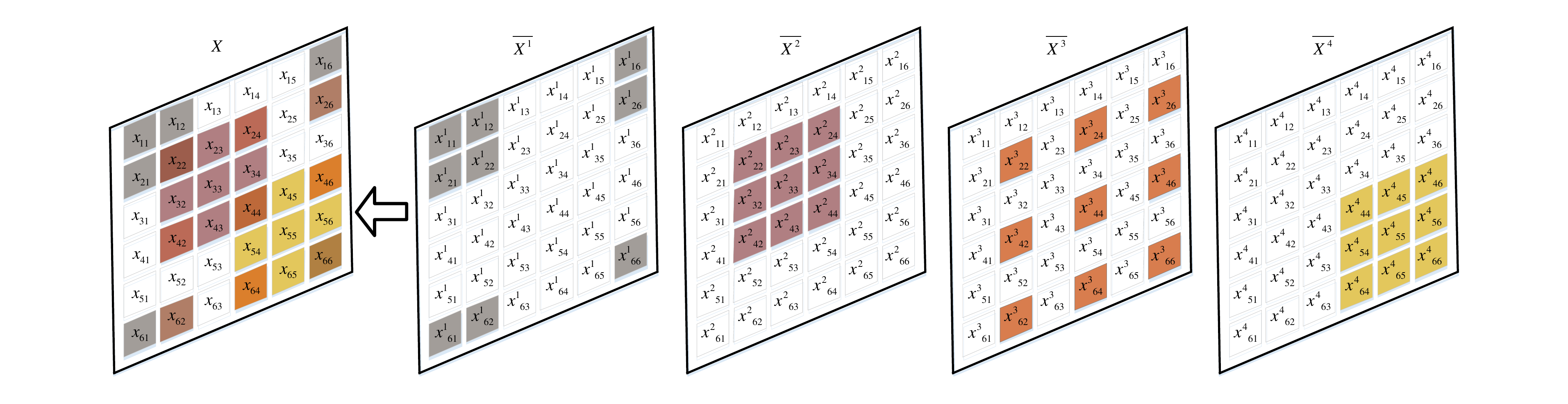}
		\caption{Matrix decomposition of global variable $X$ over four agents. }\label{graph_cut}
	\end{figure*}
\end{example}
%

\begin{remark}
	The centralized version of problem \eqref{sdp_pro} is 
	\begin{align}\label{cen-sdp}
	\min_{X\in\mathbb{S}_{+}^n} \,\, \langle M, X\rangle,
	\quad {\rm s. \, t.}\, \,  X_{i,i}=1,\quad i\in\{1,\ldots,n\},
	\end{align}
	which is a special case of generic semi-definite programming. It appears as a convex relaxation to many problems, such as the maximum cut (MAXCUT) problems\upcite{1995Improved}, community detection\upcite{2016community} and image segmentation\upcite{mincut_graph}. In practical problems such as roadmaps or social networks, the dimension $n$ may be several millions or even billions, which makes centralized computation hard. Hence, the development of distributed algorithms for \eqref{sdp_pro} is of great importance. 
\end{remark}

\begin{remark} 
	There are two scenarios in which the problem (\ref{sdp_pro}) arises. In the first scenario, an arbitrary sparse SDP problem in the standard centralized form is converted into a distributed SDP with multiple positive semi-definite matrices $X^{\tilde{i}}$ by the idea of chordal decomposition of positive semi-definite cones in \cite{chor_spar}. In the second scenario, it is assumed that the SDP is associated with a multi-agent network and matches the formulation in (\ref{sdp_pro}) exactly, such as large-scale image segmentation by multiple agents in section Simulation.
\end{remark}

Since $\langle\overline{M}^{\tilde{i}},  X\rangle$ only depends on  elements with indices in $J_{\tilde{i}}\times J_{\tilde{i}}$ of variable $X\in\mathbb{S}_{+}^{n}$, define matrix $M^{\tilde{i}}\in \mathbb{S}^{|J_{\tilde{i}}|}$ ($X^{\tilde{i}}\in \mathbb{S}^{|J_{\tilde{i}}|}$) as the remaining matrix by deleting elements of $\overline{M}^{\tilde{i}}\in\mathbb{S}_{+}^{n}$ ($\overline{X}^{\tilde{i}}\in\mathbb{S}_{+}^{n}$), whose indices are not in $J_{\tilde{i}}\times J_{\tilde{i}}$. For ease of notation, we define $X^{\tilde{i}}_{\{j,k\}}\triangleq \overline{X}^{\tilde{i}}_{(j,k)}$ for $(j,k)\in J_{\tilde{i}}\times J_{\tilde{i}}$ and $\tilde{i}\in\mathcal V$. Similarly, define $M^{\tilde{i}}_{\{j,k\}}\triangleq \overline{M}^{\tilde{i}}_{(j,k)}$ for $(j,k)\in J_{\tilde{i}}\times J_{\tilde{i}}$ and $\tilde{i}\in\mathcal V$. Hence, problem \eqref{sdp_pro} is equivalent to
\begin{subequations}\label{sdp_pro2}
	\begin{align}
	\min_{X^{\tilde{i}}\in \mathbb{S}_+^{|J_{\tilde{i}}|},\,\tilde{i}\in\mathcal V} \,\, & \sum_{\tilde{i}=1}^m  f_{\tilde{i}}(X^{\tilde{i}}), \quad   f_{\tilde{i}}(X^{\tilde{i}})= \langle {M}^{\tilde{i}}, X^{\tilde{i}}\rangle, \label{sdp_func}\\
	{\rm s. \, t.}\, \, & X^{\tilde{i}}_{\{j,j\}}=1, \ j \in J_{\tilde{i}}, \label{sdp_constr}\\
	& X^{\tilde{i}}_{\{l,k\}}= X^{\tilde{j}}_{\{l,k\}}, \ \forall l,k \in J_{\tilde{i}}\cap J_{\tilde{j}},\ (\tilde{i},\tilde{j})\in \mathcal{E},
	\end{align}
\end{subequations}
where agent $\tilde{i}\in\mathcal V$ knows ${M}^{\tilde{i}}$ and computes a positive semi-definite matrix variable $X^{\tilde{i}}\in \mathbb{S}_{+}^{|J_{\tilde{i}}|}$, $(\tilde{i},\tilde{j})\in \mathcal{E}$ specifies an overlap between the local variables $X^{\tilde{i}}$ and $X^{\tilde{j}}$ of agents $\tilde{i}$ and $\tilde{j}$.

\par In most existing distributed works\upcite{ADMM_CDC,fast_dis} for (\ref{sdp_pro}), the updating of local variable $X^{\tilde{i}}$ often involves a projection operator to the positive semi-definite cone. If the dimension of local $X^{\tilde{i}}$ of large-scale SDP is large, the projection operator is difficult and time-consuming. Hence, based on prior works on the low-rank property of matrix variables, we further reduce the computational and storage burden by representing $X\in\mathbb R^{n\times n}$ by $V'V$ with $V=[v_1,\cdots, v_n]\in\mathbb R^{p\times n}$ to avoid the projection operator. It is well-known that the rank of an optimal solution is at most $\lceil\sqrt{2n}\rceil$ (see \cite{SDP_rank}). Let $v^{\tilde{i}}_s$ be the estimate of $v_s$ by agent $\tilde{i}$ for $s\in J_{\tilde{i}}$ and $\tilde{i}\in\mathcal V$. Without causing confusions, we define
$V^{\tilde{i}}=[v^{\tilde{i}}_{s_1},\cdots,v^{\tilde{i}}_{s_{|J_{\tilde{i}}|}}]\in\mathbb R^{p\times |J_{\tilde{i}}|}$, where  $\{s_1,\ldots,s_{|J_{\tilde{i}}|}\}= J_{\tilde{i}}$, $s_1<\cdots<s_{|J_{\tilde{i}}|}$, and $\tilde{i}\in\mathcal V$.
Hence, the local variable $X^{\tilde{i}}$ is replaced by $X^{\tilde{i}}=V^{\tilde{i}'}V^{\tilde{i}}$, where $V^{\tilde{i}}\in \mathbb{R}^{p\times \left|J_{\tilde{i}}\right|}$, $p>\sqrt{2n}$. Then the semi-definite programming (\ref{sdp_pro2}) is rewritten as the following non-convex optimization problem on unit spheres:
\begin{subequations}\label{non-convex_pro}
	\begin{align}
	\min_{V^{\tilde i},\,\tilde{i}\in\mathcal V} \quad & \sum_{\tilde{i}=1}^m \tilde f_{\tilde{i}}(V^{\tilde{i}}), \ \tilde  f_{\tilde{i}}(V^{\tilde{i}})=\langle M^{\tilde{i}} ,V^{\tilde{i}'}V^{\tilde{i}}\rangle \label{obj}\\
	{\rm s.} \ {\rm t.} \quad  &\|v_j^{\tilde{i}}\|=1, \quad \forall \tilde{i}\in\mathcal V, \ {j} \in J_{\tilde{i}}\label{norm_con}\\
	&v_{{s}}^{\tilde{i}}=v_{{s}}^{\tilde{j}}, \quad \forall {s} \in J_{\tilde{i}}\cap J_{\tilde{j}},\ (\tilde{i},\tilde{j})\in \mathcal{E},\label{share_con}
	\end{align}
\end{subequations}
where $V^{\tilde{i}}\in \mathbb{R}^{p \times \left|J_{\tilde{i}}\right|}$ is local variable and $M^{\tilde{i}}\in\mathbb{S}^{\left|J_{\tilde{i}}\right|}$ is local coefficient matrix known by agent $\tilde{i}\in\mathcal V$.
\par The following assumption is needed.
\begin{assumption}\label{pro_assump}
	\begin{itemize}
		\item[(1)] Graph $\mathcal G$ is undirected and connected.
		\item[(2)] Each element of global coefficient matrix $M=\sum_{\tilde{i}=1}^m \overline{M}^{\tilde{i}}$ is non-negative and  diagonal elements of $M$ are zero.
		\item[(3)] The integer $p$  satisfies $p>\sqrt{2n}$.
	\end{itemize}
\end{assumption}
\begin{remark}
	Assumption \ref{pro_assump} (1) is general in distributed optimization.
	Since the norm of column variables $v_i$ is fixed as one,  $M_{(i,i)}=0$ does not affect the solution of optimization problem. In addition,  because elements of coefficient matrix  in MAXCUT problems and community detection are non-negative, 
	Assumption \ref{pro_assump} (2) is  practical for   problem \eqref{sdp_pro}.  Assumption \ref{pro_assump} (3) is a sufficient condition that optimal solutions for $V^{\tilde{i}}$'s recover optimal solutions for $X^{\tilde{i}}$'s.
\end{remark}

\subsection{Distributed synchronous optimization algorithm}\label{pro_trans}
In this subsection, we propose a distributed synchronous algorithm for the transformed distributed non-convex optimization problem (\ref{non-convex_pro}).

To present the algorithm, we need additional definitions and notations. We define ``children" and ``parents" in graph $\mathcal G$. For the problem (\ref{non-convex_pro}), there is an edge between agents $\tilde{i}$ and $\tilde{j}$ if $J_{\tilde{i}}\cap J_{\tilde{j}} \neq \emptyset$. Without loss of generality,  if indices $\tilde{j}<\tilde{i}$ and $J_{\tilde{i}}\cap J_{\tilde{j}} \neq \emptyset$, the agent $\tilde{j}$ is called the {\em parent} of agent $\tilde{i}$, denoted by ${\rm par}(\tilde{i})$, and the set of {\em children} of agent $\tilde{i}$ is denoted by ${\rm ch}(\tilde{i})$. The message passed from agent $\tilde{r}\in {\rm ch}(\tilde{i})$ to agent $\tilde{i} \in\mathcal V$ is denoted by $\varpi^{\tilde{r},\tilde{i}}$. Define the {\em coupling set} of indices between column variables of agent $\tilde{i}$ and its parent ${\rm par}(\tilde{i})$ as $\mathcal{S}_{\tilde{i},{\rm par}(\tilde{i})}={{J}_{\tilde{i}}}\cap J_{{\rm par}(\tilde{i})}$ and the {\em uncoupling set} of indices as $\mathcal{R}_{\tilde{i},{\rm par}(\tilde{i})}=J_{\tilde{i}}\backslash \mathcal{S}_{\tilde{i},{\rm par}(\tilde{i})}$. In addition,  column variables $v_j^{\tilde{i}}$s with $j\in \mathcal{S}_{\tilde{i},{\rm par}(\tilde{i})}$ are called {\em coupling variables} of agent $\tilde{i}$ and its parent ${\rm par}(\tilde{i})$, and  column variables $v_j^{\tilde{i}}$'s with $j\in \mathcal{R}_{\tilde{i},{\rm par}(\tilde{i})}$ are called {\em uncoupling variables} of agent $\tilde{i}$ and its parent ${\rm par}(\tilde{i})$. Accordingly, we have the similar notations $\mathcal{S}_{{\rm ch}(\tilde{i}),\tilde{i}}$ and $\mathcal{R}_{{\rm ch}(\tilde{i}),\tilde{i}}$.

\begin{example}
	Fig. \ref{parent_child_graph} gives  a network  to show the previous definitions of $\mathcal{S}_{\cdot,\cdot}$ and $\mathcal{R}_{\cdot,\cdot}$. Noted that $\mathcal{S}_{{\rm ch}(\tilde{i}),\tilde{i}}$ is an empty set for agent $\tilde{i}$ with no children, such as agents $\tilde{3},\tilde{4},\tilde{5}$. Accordingly, the message $\varpi^{\tilde{r},\tilde{i}}, \ \tilde{r}\in {\rm ch}(\tilde{i})$ passed from children agents is zero for $\tilde{i}=\tilde{3},\tilde{4},\tilde{5}$. The similar case $\mathcal{S}_{\tilde{i},{\rm par}(\tilde{i})}=\emptyset$ holds for the root agent, which have no parents. 
	\begin{figure}[h]
		\centering
		\includegraphics[width=8cm]{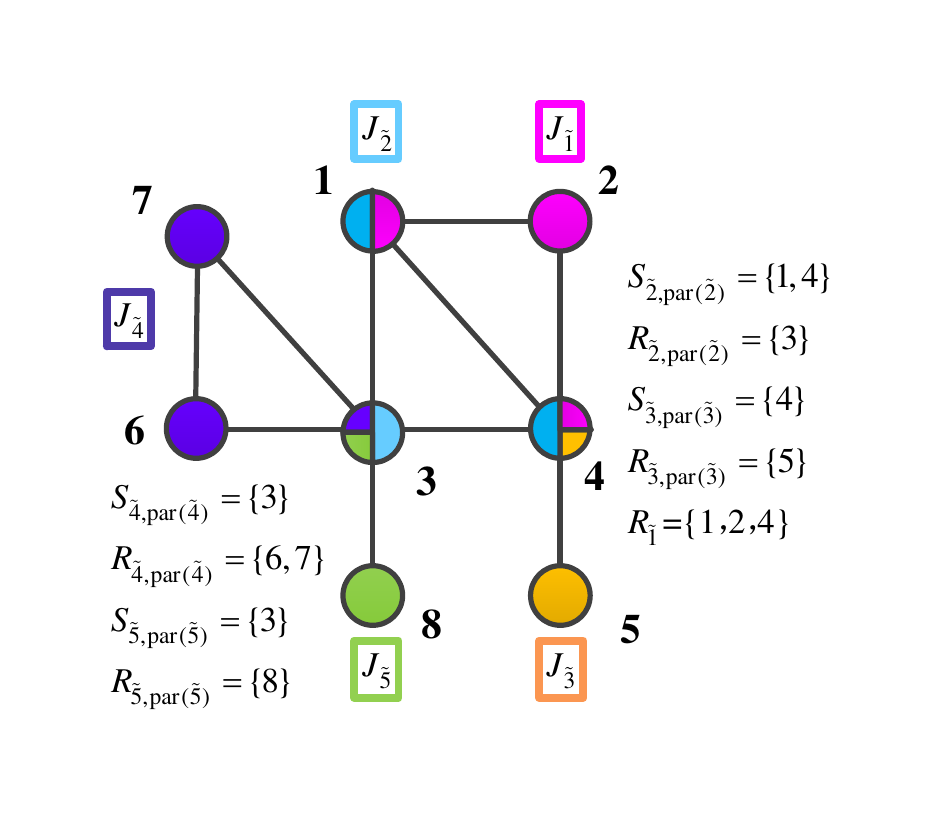}\\
		\caption{The coupling and uncoupling variables over a sparse graph. The  subscript elements $\tilde{(\cdot)}$ of set $\mathcal{S}$ or $\mathcal{R}$ denote different agents and the elements $(\cdot)$ in the brace denote the indices of coupling variables among different agents. } \label{parent_child_graph}
		\centering
	\end{figure}
\end{example}

\par  For each agent $\tilde{i}\in\mathcal V$ and $j\in \mathcal{R}_{\tilde{i},{\rm par}(\tilde{i})}$, define $$p_j^{\tilde{i}}(t)\triangleq\sum_{l<j,l\in J_{\tilde{i}}} M_{\{j,l\}}^{\tilde{i}} {v}_l^{\tilde{i}}(t+1)+\sum_{l>j,l\in J_{\tilde{i}}} M_{\{j,l\}}^{\tilde{i}} {v}_l^{\tilde{i}}(t).$$ 
Define  the step-size $\theta_j^{\tilde{i}}$ as
\begin{align}\label{theta_de}
\theta_j^{\tilde{i}} \in(0,\frac{1}{\sum_{\tilde{s}\in (\tilde{i}, {\rm ch}(\tilde{i}))}\|M_{\{j,:\}}^{\tilde{s}}\|_1}),\ j\in \mathcal{R}_{\tilde{i},{\rm par}(\tilde{i})}
\end{align}
where $M_{\{j,:\}}^{\tilde{s}}$ denotes the $j$th row of matrix $\overline{M}^{\tilde{s}}$.
\par The massage passed from agent $\tilde{r}\in {\rm ch}(\tilde{i})$ to its parent  $\tilde{i}\in \mathcal{V}$ is defined as
\begin{align}\label{distri_mes_syn}
\varpi_j^{\tilde{r},\tilde{i}}(t+1) \triangleq  \sum_{l\in J_{\tilde{r}}} M_{\{j,l\}}^{\tilde{r}} v_l^{\tilde{r}}(t+1) \quad \forall j \in \mathcal{S}_{\tilde{r},\tilde{i}}.
\end{align}
and $\varpi_j^{\tilde{r},\tilde{i}}(t+1)=0$ for $j\in \mathcal{R}_{\tilde{i},{\rm par}(\tilde{i})}\backslash \mathcal{S}_{\tilde{r},\tilde{i}}$.

For $\tilde i\in\mathcal V$ and $j\in \mathcal{R}_{\tilde{i},{\rm par}(\tilde{i})}$, the variable $v_j^{\tilde{i}}(t+1)$  is updated as
\begin{align}\label{vupdatesm}
&v_j^{\tilde{i}}(t+1)={\rm normal}\Big(v_j^{\tilde{i}}(t)-\theta_j^{\tilde{i}} \big[p_j^{\tilde{i}}(t)+\sum_{\tilde{r}\in {\rm ch}(\tilde{i})}\varpi_j^{\tilde{r},\tilde{i}}(t+1)\big]\Big).
\end{align}

The message sent from agent $\tilde i\in\mathcal V$ to its parent ${\rm par}(\tilde{i})$ is
\begin{align}\label{message2par}
\varpi_j^{\tilde{i},{\rm par}(\tilde{i})} =\sum_{l\in J_{\tilde{i}}} M_{\{j,l\}}^{\tilde{i}} v_l^{\tilde{i}}(t+1)+\sum_{\tilde{r} \in {\rm ch}(\tilde{i})}\varpi_j^{\tilde{r},\tilde{i}}(t+1),
\end{align}		
where $j\in \mathcal{S}_{\tilde{i},{\rm par}(\tilde{i})}$.

Define $\mathbf V = (V^{1},\cdots,V^{m})$. The distributed synchronous algorithm is given in Algorithm \ref{algo_sum}.

\begin{algorithm}[H]
	\caption{Distributed Synchronous Algorithm (DSA)}
	\label{algo_sum}
	\begin{algorithmic}[1]  
		\State  \textbf{Initialization:}  Initialize  $v_j^{\tilde{i}}=v_0\in\mathbb R^p$ such that $\|v_0\|=1$ for all  $\tilde{i}\in\mathcal V$ and $j\in J_{\tilde{i}}$.
		\While{the stopping criteria is not satisfied}
		\For{$\tilde{i}=m$ to $\tilde{i}=1$}
		\For {$j\in \mathcal{R}_{\tilde{i},{\rm par}(\tilde{i})}$}
		\If {$j\in \mathcal{S}_{{\rm ch}(\tilde{i}),\tilde{i}}$}
		\State agent $\tilde{i}$ receives $\varpi_j^{\tilde{r},\tilde{i}}(t+1)$ (computed by \eqref{distri_mes_syn}) from its child $\tilde{r}\in {\rm ch}(\tilde{i})$.
		\EndIf
		\State   agent $\tilde{i}$ updates variable $v^{\tilde{i}}_j(t+1)$ following \eqref{vupdatesm}.
		\If{$j\in \mathcal{S}_{{\rm ch}(\tilde{i}),\tilde{i}}$}
		\State  agent $\tilde{i}$ sends $v_j^{\tilde{i}}(t+1)$ to its child $\tilde{r}\in{\rm ch}(\tilde{i})$, i.e., 
		\begin{align}\label{syn_sup}
			v_j^{\tilde{r}}(t+1)=v_j^{\tilde{i}}(t+1).
		\end{align}
		\EndIf
		\EndFor
		\For {$j\in \mathcal{S}_{\tilde{i},{\rm par}(\tilde{i})} $}
		\If {$j\in \mathcal{S}_{{\rm ch}(\tilde{i}),\tilde{i}}$}
		\State agent $\tilde{i}$ receives information $\varpi_j^{\tilde{r},\tilde{i}}$ from its child $\tilde{r}\in {\rm ch}(\tilde{i})$.
		\EndIf
		
		\State agent $\tilde{i}$ sends message $\varpi_j^{\tilde{i},{\rm par}(\tilde{i})}$ (computed by \eqref{message2par}) to its parent ${\rm par}(\tilde{i})$.
		
		\EndFor
		\EndFor
		\State $t\leftarrow t+1$.
		\EndWhile
	\end{algorithmic}
\end{algorithm}

The convergence performance of Algorithm 1 is provided in the following theorem, whose proof is given in the next section.
\begin{theorem}\label{syn_theo}
	Let Assumption \ref{pro_assump} hold and $\{\mathbf V(t)\}$ be a sequence generated by Algorithm 1. Then $\mathbf V(t)$ converges to global optimal solutions to (\ref{non-convex_pro}) almost surely under random initialization as $t\rightarrow\infty$. 
\end{theorem}
\begin{remark}
	Note that there may be several parents for each agent in practice. In this paper, we only consider the case that each agent owns one parent for convenience of analysis. However, the algorithm can be  easily extended to cases where agents have multiple parents. We  provide examples in simulations to demonstrate the cases of tree-structured graphs of which one agent has a unique parent and  graphs of which one agent may have multiple parents.	
\end{remark}
\begin{remark}\label{syn_mark}
	Compared with the existing centralized works \cite{wang2017mixing,saddle_escape}, the proposed algorithm decentralizes the storage space and computational burden of large-scale SDP over different agents at the cost of network communication. In addition, the proposed algorithm is applicable for the scenario where global information is located on geographically separated agents such that centralized algorithms can not handle.
\end{remark}

\subsection{Distributed asynchronous optimization algorithm}
\par The synchronous algorithm given by Algorithm \ref{algo_sum} needs a global clock and the updating rate of variables is limited by the slowest agent. Whereas, asynchronous algorithms update variables by local clocks and allow communication time-delays, then the variables updating of one agent will not be limited by other agents. Hence, in this subsection, we provide a distributed asynchronous algorithm for solving problem (\ref{non-convex_pro}).
\par Let $T^{\tilde{i}}$ be the set of times at which agent $\tilde{i}$ updates local variable $V^{\tilde{i}}$. Noted that agent $\tilde{i}$ may not have access to the most recent value of other agents' variables.
Then, we define that the updating time of variables $v_j$ at parent or child node, $\tau_{j}^{\tilde{i}}(t)$, satisfies
\begin{align}\label{taulim}
0\leq \tau_{j}^{\tilde{i}}(t) \leq t, \quad \forall t\in T^{\tilde{i}},
\end{align}
where $j\in J_{{\rm ch}(\tilde{i})}\cup J_{{\rm par}(\tilde{i})}\cup J_{\tilde{i}}$. If $j\in \mathcal{R}_{\tilde{i},{\rm par}(\tilde{i})}$, then $\tau_{j}^{\tilde{i}}(t)=t$. 
The difference $(t-\tau_{j}^{\tilde{i}}(t))$ between the current time $t$ and $\tau_{j}^{\tilde{i}}(t)$ is viewed as a form of communication delay between agent $\tilde{i}$ and its parents or children.
\par For the communication delay between agents, we assume that the following condition holds.
{\begin{assumption}\label{total_asyn}
		(Partial Asynchronism) There exists a positive integer $B$ such that:
		\begin{itemize}
			\item For each agent $\tilde{i}$ and $t\geq 0$, at least one element of the set $\{t,t+1,\cdots,t+B-1\}$ belongs to $T^{\tilde{i}}$.
			\item There holds $${\rm max}\{0,t-B+1\}\leq \tau_{{j}}^{\tilde{i}}(t)\leq t,$$
			for all agent $\tilde{i}\in \mathcal{V}$ and ${j}\in \{1,\cdots,n\}$ and all $t\geq 0$.
		\end{itemize}
\end{assumption}}

\par In the distributed asynchronous Algorithm \ref{algo_asyn}, local variable $V^{\tilde{i}}$ is updated by the time-delayed messages communicated from neighbors and local information. For each agent $\tilde{i}$, define $$p_j^{\tilde{i}}(t)= \sum_{l\in J_{\tilde{i}}} M_{\{j,l\}}^{\tilde{i}} {v}_l^{\tilde{i}}(t),$$ and the step-size $\theta_j^{\tilde{i}}$ satisfies
\begin{align}\label{asyn_step}
\theta_j^{\tilde{i}}\in (0,\frac{1}{(1+B+nB)L}),
\end{align}
where $L= \max_{\tilde{i}\in \mathcal{V}}\{l_1,\cdots,l_m\}$, $l_{\tilde{i}}$ is the best Lipschitz constant of $\nabla \tilde{f}_{\tilde{i}}(V^{\tilde{i}})$ ($\tilde{f}_{\tilde{i}}(V^{\tilde{i}})$ was defined in \eqref{non-convex_pro}) and $B$ is defined in Assumption \ref{total_asyn}.
\par The massage passed from child $\tilde{r}\in {\rm ch}(\tilde{i})$ to parent  $\tilde{i}$ is defined as
\begin{align*}
\varpi_j^{\tilde{r},\tilde{i}}(\tau^{\tilde{i}}(t)) = \sum_{l\in J_{\tilde{r}}} M_{\{j,l\}}^{\tilde{r}} v_l^{\tilde{r}}(\tau_l^{\tilde{i}}(t)) \quad \forall j \in \mathcal{S}_{\tilde{r},\tilde{i}},
\end{align*}
and for variables with indices $j\in \mathcal{R}_{\tilde{i},{\rm par}(\tilde{i})}\backslash \mathcal{S}_{\tilde{r},\tilde{i}}$,
$\varpi_j^{\tilde{r},\tilde{i}}(\tau^{\tilde{i}}(t))=0$.
\par If $t\in T^{\tilde{i}}$, the variable $v_j^{\tilde{i}}(t+1)$ is updated as
\begin{subequations}\label{v_updateasm}
\begin{align}
&v_j^{\tilde{i}}(t+1)=\!{\rm normal}\Big(v_j^{\tilde{i}}(t)-\theta_j^{\tilde{i}} \big[p_j^{\tilde{i}}(t)+\!\sum_{\tilde{r}\in {\rm ch}(\tilde{i})}\!\varpi_j^{\tilde{r},\tilde{i}}(\tau^{\tilde{i}}(t))\big]\Big),\notag\\ 
& \qquad \qquad \qquad  j\in \mathcal{R}_{\tilde{i},{\rm par}(\tilde{i})}\label{asyn_r},\\
& v_j^{\tilde{i}}(t+1)=v_j^{{\rm par}(\tilde{i})}(\tau_{j}^{\tilde{i}}(t)),\quad j\in \mathcal{S}_{\tilde{i},{\rm par}(\tilde{i})}\label{asyn_s}.
\end{align}
\end{subequations}

For times $t \notin T^{\tilde{i}}$, the variable $v_j^{\tilde{i}}$ is unchanged,

$$v_j^{\tilde{i}}(t+1)=v_j^{\tilde{i}}(t), \ \forall j\in J_{\tilde{i}}.
$$
\par The message sent from agent $\tilde{i}\in\mathcal{V}$ to its parent ${\rm par}(\tilde{i})$ is
\begin{align}\label{mess_asyn}
\varpi_j^{\tilde{i},{\rm par}(\tilde{i})} \!=\! \sum_{l\in J_{\tilde{i}}} M_{\{j,l\}}^{\tilde{i}} v_{l}^{\tilde{i}}(t+1)+\sum_{\tilde{r} \in {\rm ch}(\tilde{i})}\varpi_j^{\tilde{r},\tilde{i}}(\tau^{\tilde{i}}(t)),\ j\in J_{\tilde{i}}.
\end{align}	
\begin{algorithm}[H]
	\caption{Distributed Asynchronous Algorithm (DAA) - from the view of agent $\tilde{i}$}
	\label{algo_asyn}
	\begin{algorithmic}[1]  
		\State \textbf{Initialization:} Initialize  $v_j^{\tilde{i}}=v_0\in\mathbb R^p$ such that $\|v_0\|=1$ for all  $\tilde{i}\in\mathcal V$ and $j\in J_{\tilde{i}}$.
		\While{the stopping criteria is not satisfied}
		\State for each agent $\tilde{i}$, keep receiving information $\varpi_j^{\tilde{r},\tilde{i}}(\tau^{\tilde{i}}(t))$ from children and receiving information  $v_j^{{\rm par}(\tilde{i})}(\tau_{j}^{\tilde{i}}(t))$ from parent.
		\If{$t \in T^{\tilde{i}}$}
		\For{$j\in J_{\tilde{i}}$}		
		\State agent $\tilde{i}$ updates variable $v_j^{\tilde{i}}(t+1)$ following (\ref{v_updateasm}).
		\If{$j\in \mathcal{S}_{\tilde{i},{\rm par}(\tilde{i})}$}
		\State agent $\tilde{i}$ sends message $\varpi_j^{\tilde{i},{\rm par}(\tilde{i})}$ (computed by (\ref{mess_asyn})) to its parent,
		\EndIf
		\If{$j\in \mathcal{S}_{\tilde{l},\tilde{i}}$, $\tilde{l}\in {\rm ch}(\tilde{i})$}
		\State agent $\tilde{i}$ sends local variable $v_j^{\tilde{i}}(t+1)$ to each child $\tilde{l}$ that has coupling variable $v_j^{\tilde{i}},j\in \mathcal{S}_{\tilde{l},\tilde{i}}$.
		\EndIf
		\EndFor
		\EndIf
		\State $t\leftarrow t+1$.
		\EndWhile
	\end{algorithmic}
\end{algorithm}


\par Before providing the convergence performance of DAA, we need one additional assumption, which is vital in the transformation of proposed algorithm.
\begin{assumption}\label{asyn_assump}
	Each element of matrix $M$ is only accessible to one agent, which implies that for each agent $\tilde{i}\in \mathcal{V}$, $M^{\tilde{i}}_{\{l,k\}}=M_{(l,k)}$ holds for $(l,k)\in J_{\tilde{i}}\times J_{\tilde{i}}$.
\end{assumption}
\par Next, the convergence performance of DAA is provided in the following theorem, whose analysis is shown in the section \ref{ASM_sec}.
\begin{theorem}\label{asyn_theo}
	Under Assumptions \ref{pro_assump}-\ref{asyn_assump}, the sequence $\{\mathbf V(t)\}$ generated by DAA converges to critical points to (\ref{non-convex_pro}) as $t\rightarrow\infty$. 
\end{theorem}
\begin{remark}
	In the implementation of DAA, each agent receives communication information from its neighbors and stores it in local buffer. The local received data may be out-of-date due to time-delays. Each agent $\tilde{i}\in \mathcal{V}$ updates local variables using data in local buffer at the time $t\in T^{\tilde{i}}$ and does not have to wait for the point when other local communicating messages become available. It allows some agents to compute faster and execute more iterations than others.
	
\end{remark}
\section{Theoretical analysis}\label{proof_sec}
In this section, we present theoretical proofs for the convergence properties of proposed distributed synchronous and asynchronous algorithms, respectively.

\subsection{Convergence analysis for synchronous algorithm}
\par Firstly, we  develop a compact form of the proposed synchronous algorithm containing all agent's updates.

For any $j\in\{1,\ldots,n\}$, let $v_j \triangleq v_j^{\tilde i}$, where  $j\in \mathcal{R}_{\tilde{i},{\rm par}(\tilde{i})}$ and $v_j^{\tilde i}$ is a column vector of the variable of agent $\tilde i\in\mathcal V$. Then we define the global variable  
\begin{align}\label{V_def}
V=\left[v_1,\cdots,v_n\right]\in \mathbb{R}^{p \times n}.
\end{align}

\begin{remark}
	 For $j\in \mathcal{S}_{\tilde{i},{\rm par}(\tilde{i})}$, by the updating design \eqref{syn_sup}, there is a parent of $\tilde{i}$, $\tilde{j}\in {\rm par}(\tilde{i})$, such that $j\in \mathcal{R}_{\tilde{j},{\rm par}(\tilde{j})}$ and $v_j^{\tilde{i}}(t)=v_j^{\tilde{j}}(t)$. Hence, in the following analysis, for convenience, we consider the global variable $V$ instead of $\mathbf{V}$.
\end{remark}
 \par Without loss of generality,  we make the following assumptions in the analysis.
 \begin{assumption}\label{SR_indice}
 	Take any  $\tilde{i}\in\mathcal V$. Indices
 	$j\in \mathcal{R}_{{\rm ch}(\tilde{i}),\tilde{i}}$, $l\in \mathcal{R}_{\tilde{i},{\rm par}(\tilde{i})}$ and $p\in\mathcal{S}_{\tilde{i},{\rm par}(\tilde{i})}$ satisfy that $p>l>j$.
 \end{assumption}
 
 Notice that for each column vector $v_j$ with index $j\in \mathcal{S}_{\tilde{i},{\rm par}(\tilde{i})}$, there must exist one agent $\tilde{j}$ such that $j\in \mathcal{R}_{\tilde{j},{\rm par}(\tilde{j})}$.  
\begin{assumption}\label{twoshare_assu}
	 We assume that there is no shared variables between agents ${\rm ch}(\tilde{i})$ and ${\rm par}(\tilde{i})$ for any $\tilde{i}\in\mathcal V$. That is,
 $\mathcal{S}_{{\rm ch}(\tilde{i}),\tilde{i}}\subset \mathcal{R}_{\tilde{i},{\rm par}(\tilde{i})}$ for all $\tilde{i}\in\mathcal V$.
\end{assumption}
\par Note that Assumptions \ref{SR_indice} and \ref{twoshare_assu} are not needed in the proposed synchronous algorithm, which is developed for the convenience of proof. 

\par Define $\Theta\in\mathbb{R}^{n\times n}$ as a diagonal matrix with diagonal elements,
\begin{align}\label{dia_theta}
\Theta_{(j,j)}=\theta_j^{\tilde{i}} \quad {\rm if} \ j\in \mathcal{R}_{\tilde{i},{\rm par}(\tilde{i})}.
\end{align}

Then, under the above assumptions, we obtain the following result.
\begin{lemma}\label{cmp_lem}
	Under Assumptions \ref{SR_indice} and \ref{twoshare_assu}, the distributed synchronous algorithm is equivalent to
	\begin{align}\label{distri_mapping}
	v_j(t+1)=&{\rm normal}\big(v_j(t)-\Theta_{(j,j)} \mathbf{g}_j(t)\big), \ \forall \ j\in \{1,\cdots,n\},
	\end{align}
	where $\mathbf{g}_j(t)=\sum_{l\in  \{1,\cdots,j-1\}} M_{(j,l)} v_l(t+1)+\sum_{l\in \{j+1,\cdots,n\}} M_{(j,l)} v_l(t)$.
\end{lemma}
\par\textbf{Proof:}
For variables $v_j^{\tilde{i}}$ of agent $\tilde{i}$ such that $j\in \mathcal{R}_{\tilde{i},{\rm par}(\tilde{i})}$, we substitute (\ref{distri_mes_syn}) to the updating (\ref{vupdatesm}) and get
\begin{align}\label{vr_up}
v_j^{\tilde{i}}(t+1)
={\rm normal}\Big(v_j^{\tilde{i}}(t)-\theta_j^{\tilde{i}} \mathbf{g}^{\tilde{i}}_j(t)\Big),\ j\in \mathcal{R}_{\tilde{i},{\rm par}(\tilde{i})}
\end{align}
where 
\begin{align}\label{PHI}
\mathbf{g}^{\tilde{i}}_j(t)=&\sum_{l\in  \{1,\cdots,j-1\}} M^{\tilde{i}}_{\{j,l\}} {v}_l^{\tilde{i}}(t+1)\notag\\
&+\sum_{l\in \{j+1,\cdots,n\}} M^{\tilde{i}}_{\{j,l\}} {v}_l^{\tilde{i}}(t)\notag\\
&+\sum_{\tilde{r}\in {\rm ch}(\tilde{i})}\sum_{l\in J_{\tilde{r}}} M_{\{j,l\}}^{\tilde{r}} v_{l}^{\tilde{r}}(t+1).
\end{align}
\par The condition $j\in \mathcal{R}_{\tilde{i},{\rm par}(\tilde{i})}$ holds in the whole process of proof. For convenience, in the following analysis, we omit this condition.
\par {Since} $M_{\{j,l\}}^{\tilde{r}}=\overline{M}_{\{j,l\}}^{\tilde{r}}$, for $(j,l)\in J_{\tilde{r}}\times J_{\tilde{r}}$ and $\tilde{r}\in \mathcal{V}$, we have
\begin{align}\label{msum}
&\sum_{\tilde{r}\in \{\tilde{i},{\rm ch}(\tilde{i})\}} {M}^{\tilde{r}}_{\{j,l\}}=\sum_{\tilde{r}=1}^m \overline{M}^{\tilde{r}}_{\{j,l\}}=M_{(j,l)},\\
& j\in \mathcal{R}_{\tilde{i},{\rm par}(\tilde{i})}\subset  \cup_{\tilde{r}\in \{\tilde{i},{\rm ch}(\tilde{i})\}} J_{\tilde{r}}, l\in \cup_{\tilde{r}\in \{\tilde{i},{\rm ch}(\tilde{i})\}}J_{\tilde{r}},\notag
\end{align}
where the first equality holds because for agent $\tilde{j}\in \{1,\cdots,m\}\backslash\{\tilde{i},{\rm ch}(\tilde{i})\}$, $M_{\{j,l\}}^{\tilde{j}}=0$, where $j\in \mathcal{R}_{\tilde{i},{\rm par}(\tilde{i})}\subset  \cup_{\tilde{r}\in \{\tilde{i},{\rm ch}(\tilde{i})\}} J_{\tilde{r}}, l\in \cup_{\tilde{r}\in \{\tilde{i},{\rm ch}(\tilde{i})\}}J_{\tilde{r}}$.
\par By \eqref{msum}, the term $\sum_{\tilde{r}\in {\rm ch}(\tilde{i})}\sum_{l\in J_{\tilde{r}}} M_{\{j,l\}}^{\tilde{r}} v_{l}^{\tilde{r}}(t+1)$ in \eqref{PHI} satisfies
\begin{align}\label{lastterm_trans}
&\sum_{\tilde{r}\in {\rm ch}(\tilde{i})}\sum_{l\in J_{\tilde{r}}} M_{\{j,l\}}^{\tilde{r}} v_{l}^{\tilde{r}}(t+1) \notag\\
=&\sum_{\tilde{r}\in {\rm ch}(\tilde{i})}\!\sum_{l\in \mathcal{S}_{\tilde{r},\tilde{i}}}\! M_{\{j,l\}}^{\tilde{r}} v_{l}^{\tilde{r}}(t+1)\!+\!\sum_{\tilde{r}\in {\rm ch}(\tilde{i})}\!\sum_{l\in \mathcal{R}_{\tilde{r},\tilde{i}}}\! M_{\{j,l\}}^{\tilde{r}} v_{l}^{\tilde{r}}(t+1)\notag\\
=&\sum_{\tilde{r}\in {\rm ch}(\tilde{i})}\!\sum_{l\in \mathcal{S}_{\tilde{r},\tilde{i}}}\! M_{\{j,l\}}^{\tilde{r}} v_{l}^{\tilde{r}}(t+1)\!+\!\sum_{\tilde{r}\in {\rm ch}(\tilde{i})}\!\sum_{l\in \mathcal{R}_{\tilde{r},\tilde{i}}}\! M_{\{j,l\}}^{\tilde{r}} v_{l}(t+1)\notag\\
=&\sum_{\tilde{r}\in {\rm ch}(\tilde{i})}\sum_{l\in \mathcal{S}_{\tilde{r},\tilde{i}}\cap \{l<j\}} M_{\{j,l\}}^{\tilde{r}} v_{l}^{\tilde{r}}(t+1)\notag\\
&+\sum_{\tilde{r}\in {\rm ch}(\tilde{i})}\sum_{l\in \mathcal{S}_{\tilde{r},\tilde{i}}\cap \{l>j\}} M_{\{j,l\}}^{\tilde{r}} v_{l}^{\tilde{r}}(t)\notag\\
&+\sum_{l\in \mathcal{R}_{\tilde{r}\in {\rm ch}(\tilde{i}),\tilde{i}}} M_{(j,l)} v_{l}(t+1)
\end{align}
where the last equality holds because variables $v_l^{\tilde{r}}$ with indices satisfying $l>j,j\in \mathcal{R}_{\tilde{i},{\rm par}(\tilde{i})}$ are not updated by the algorithm design and $M^{\tilde{i}}_{\{j,l\}}=0$ for $l\in \mathcal{R}_{\tilde{r}{\tilde{i}}}$.
\par Then, by substituting \eqref{lastterm_trans} to $\mathbf{g}_j^{\tilde{i}}$ in \eqref{PHI}, we obtain
\begin{align}\label{Phi_com}
\mathbf{g}^{\tilde{i}}_j(t)=\zeta_j(t+1)+\iota_j
(t)
\end{align} 
where 
\begin{align*}
\zeta_j(t+1)=&\sum_{l\in  \{1,\cdots,j-1\}} M^{\tilde{i}}_{\{j,l\}} {v}_l^{\tilde{i}}(t+1)\\
&+\sum_{\tilde{r}\in {\rm ch}(\tilde{i})}\sum_{l\in \mathcal{S}_{\tilde{r},\tilde{i}}\cap \{l<j\}} M_{\{j,l\}}^{\tilde{r}} v_{l}^{\tilde{r}}(t+1)\\
&+\sum_{l\in \mathcal{R}_{\tilde{r}\in {\rm ch}(\tilde{i}),\tilde{i}}} M_{(j,l)} v_{l}(t+1),
\end{align*} 
and 
\begin{align*}
\iota_j
(t)=&\sum_{l\in \{j+1,\cdots,n\}} M^{\tilde{i}}_{\{j,l\}} {v}_l^{\tilde{i}}(t)\\
&+\sum_{\tilde{r}\in {\rm ch}(\tilde{i})}\sum_{l\in \mathcal{S}_{\tilde{r},\tilde{i}}\cap \{l>j\}} M_{\{j,l\}}^{\tilde{r}} v_{l}^{\tilde{r}}(t).
\end{align*}
\par  Then, we will discuss variables $\zeta_j(t+1)$ and $\iota_j(t)$ respectively. 
\par (1) Consider $\zeta_j$ composed of variables with indices $l$ in the set $\{l|l<j,j\in \mathcal{R}_{\tilde{i},{\rm par}(\tilde{i})}\}$. Let $MR\triangleq\sum_{l\in \mathcal{R}_{\tilde{i},{\rm par}(\tilde{i})}\cap \{l<j\}} M^{\tilde{i}}_{\{j,l\}} {v}_l^{\tilde{i}}(t+1)$ and $MS\triangleq\sum_{\tilde{r}\in {\rm ch}(\tilde{i})}\sum_{l\in \mathcal{S}_{\tilde{r},\tilde{i}}\cap \{l<j\}} M_{\{j,l\}}^{\tilde{r}} v_{l}^{\tilde{r}}(t+1)$. The sum of $MS$ and $MR$ is
\begin{align}\label{lj2}
&MS+MR\notag\\
=&\sum_{\tilde{r}\in \{\tilde{i},{\rm ch}(\tilde{i})\}}\sum_{l\in  \mathcal{R}_{\tilde{i},{\rm par}(\tilde{i})} \cap \{l<j\}} M_{\{j,l\}}^{\tilde{r}} v_{l}^{\tilde{r}}(t+1)\notag\\
=&\sum_{\tilde{r}\in \{\tilde{i},{\rm ch}(\tilde{i})\}}\sum_{l\in  \mathcal{R}_{\tilde{i},{\rm par}(\tilde{i})} \cap \{l<j\}} M^{\tilde{r}}_{\{j,l\}}v_l(t+1)\notag\\
=&\sum_{l\in  \mathcal{R}_{\tilde{i},{\rm par}(\tilde{i})} \cap \{l<j\}} M_{(j,l)}v_l(t+1).
 \end{align}
where the first equality holds because of the condition $\mathcal{S}_{{\rm ch}(\tilde{i}),\tilde{i}}\subseteq \mathcal{R}_{\tilde{i},{\rm par}(\tilde{i})}$, and the last equality holds due to the relationship (\ref{msum}).
\par By Assumption \ref{SR_indice}, any index $l$ in $\mathcal{R}_{\tilde{r}\in {\rm ch}(\tilde{i}),\tilde{i}}$ satisfies $l<j$ for $j\in \mathcal{R}_{\tilde{i},{\rm par}(\tilde{i})}$. Then, with \eqref{lj2}, 
\begin{align}\label{lleqj}
&\zeta_j(t+1)\notag \\
=&\sum_{l\in \mathcal{R}_{\tilde{r}\in {\rm ch}(\tilde{i}),\tilde{i}}}\! M_{(j,l)} v_{l}(t+1)\!+\!\sum_{l\in \mathcal{R}_{\tilde{i},{\rm par}(\tilde{i})}\cap\{l<j\}} \! M_{(j,l)}v_l(t+1)\notag\\
=&\sum_{l\in \{l<j,j\in \mathcal{R}_{\tilde{i},{\rm par}(\tilde{i})}\}} M_{(j,l)}v_l(t+1).
\end{align}
\par (2) Consider $\iota_j$ composed of variables with indices $l$ in the set $\{l|l>j,j\in \mathcal{R}_{\tilde{i},{\rm par}(\tilde{i})}\}$. The term $\sum_{l>j} M^{\tilde{i}}_{\{j,l\}} {v}_l^{\tilde{i}}(t)$ in $\iota_j(t)$ satisfies
\begin{align}\label{lgeqj1}
&\sum_{l\in \{1,\cdots,n\}\cap\{l>j\}} M^{\tilde{i}}_{\{j,l\}} {v}_l^{\tilde{i}}(t)\notag \\
=&\sum_{l\in \mathcal{R}_{\tilde{i},{\rm par}(\tilde{i})}\cap \{l>j\}} M_{\{j,l\}}^{\tilde{i}} v_{l}^{\tilde{i}}(t)+\!\sum_{l\in \mathcal{S}_{\tilde{i},{\rm par}(\tilde{i})}\cap \{l>j\}} M_{\{j,l\}}^{\tilde{i}} v_{l}^{\tilde{i}}(t).
\end{align}
\par It follows from a similar analysis in \eqref{lj2} that the sum of $\sum_{\tilde{r}\in {\rm ch}(\tilde{i})}\sum_{l\in \mathcal{S}_{\tilde{r},\tilde{i}}\cap \{l>j\}} M_{\{j,l\}}^{\tilde{r}} v_{l}^{\tilde{r}}(t)$ in $\iota_j(t)$ and $\sum_{l\in \mathcal{R}_{\tilde{i},{\rm par}(\tilde{i})}\cap \{l>j\}} M_{\{j,l\}}^{\tilde{i}} v_{l}^{\tilde{i}}(t)$ in (\ref{lgeqj1}) is
\begin{align}\label{Rlgeqj}
&\!\sum_{\tilde{r}\in {\rm ch}(\tilde{i})}\sum_{l\in \mathcal{S}_{\tilde{r},\tilde{i}}\cap \{l>j\}}\! M_{\{j,l\}}^{\tilde{r}}  v_{l}^{\tilde{r}}\!(t)+\!\sum_{l\in \mathcal{R}_{\tilde{i},{\rm par}(\tilde{i})}\cap \{l>j\}}\! M_{\{j,l\}}^{\tilde{i}} \!v_{l}^{\tilde{i}}\!(t)\notag\\
=&\sum_{l\in  \mathcal{R}_{\tilde{i},{\rm par}(\tilde{i})} \cap \{l>j\}} M_{(j,l)}v_l(t).
\end{align}
\par Then, consider the term $\sum_{l\in \mathcal{S}_{\tilde{i},{\rm par}(\tilde{i})}\cap \{l>j\}} M_{\{j,l\}}^{\tilde{i}} v_{l}^{\tilde{i}}(t)$ in (\ref{lgeqj1}).
Because all column variables are initialized as a same value and the condition $\mathcal{S}_{{\rm ch}(\tilde{i}),\tilde{i}}\subset \mathcal{R}_{\tilde{i},{\rm par}(\tilde{i})}$ in Assumption \ref{twoshare_assu}, $v_l^{\tilde{i}}(t)=v_l^{{\rm par}(\tilde{i})}(t)=v_l(t), \ \forall  l\in \mathcal{S}_{\tilde{i},{\rm par}(\tilde{i})}$, by Algorithm 1. By (\ref{msum}), we have, for $j\in \mathcal{R}_{\tilde{i},{\rm par}(\tilde{i})}$,
\begin{align}\label{Slgeqj}
&\sum_{l\in \mathcal{S}_{\tilde{i},{\rm par}(\tilde{i})}\cap \{l>j\}}M_{\{j,l\}}^{\tilde{i}}v_l^{\tilde{i}}(t)\notag\\
=&\sum_{l\in \mathcal{S}_{\tilde{i},{\rm par}(\tilde{i})}\cap \{l>j\}}M_{(j,l)}v_l(t),
\end{align}
where the equality holds because for each child $\tilde{r}\in {\rm ch}(\tilde{i})$, $M^{\tilde{r}}_{j,l}=0$, for $j\in \mathcal{R}_{\tilde{i},{\rm par}(\tilde{i})}$, $l\in \mathcal{S}_{\tilde{i},{\rm par}(\tilde{i})}$.
\par Then, with (\ref{Rlgeqj}) and (\ref{Slgeqj}), 
\begin{align}\label{lgeqj}
&\iota_j(t)\notag\\
=&\sum_{l\in  \mathcal{R}_{\tilde{i},{\rm par}(\tilde{i})} \cap \{l>j\}} M_{(j,l)}v_l(t)+\sum_{l\in \mathcal{S}_{\tilde{i},{\rm par}(\tilde{i})}\cap \{l>j\}}M_{(j,l)}v_l(t)\notag\\
=&\sum_{l\in \{j+1,\cdots,n\}} M_{(j,l)}v_l(t).
\end{align}

\par  Hence, by \eqref{Phi_com}, (\ref{lleqj}) and  (\ref{lgeqj}), the updating (\ref{vr_up}) of agent $\tilde{i}$ for any vector variable with index $j\in \mathcal{R}_{\tilde{i},{\rm par}(\tilde{i})}$ is 
\begin{align}\label{v_comp_end}
&v_j^{\tilde{i}}(t+1)={\rm normal}\big(v^{\tilde{i}}_j(t)-\Theta_{(j,j)}\mathbf{g}^{\tilde{i}}_j(t)\big),
\end{align}
where $\mathbf{g}_j^{\tilde{i}}(t)=\sum_{l\in  \{l<j,j\in \mathcal{R}_{\tilde{i},{\rm par}(\tilde{i})}\}} M_{(j,l)} v_l(t+1)+\sum_{l\in \{j+1,\cdots,n\}\cap \{j\in \mathcal{R}_{\tilde{i},{\rm par}(\tilde{i})} \}} M_{(j,l)} v_l(t)$ and $\Theta_{(j,j)}=\theta_j^{\tilde{i}}$, which is defined in (\ref{theta_de}). Since \eqref{v_comp_end} holds for each agent $\tilde{i}$, we obtain the desire result \eqref{distri_mapping} with \eqref{V_def}. 
%
$\hfill\blacksquare$
\par Next, we will discuss the convergence properties of proposed distributed synchronous algorithm. Because the optimization problem (\ref{non-convex_pro}) is a non-convex minimization problem, there may exist several local minima and saddle points, which are regarded as major obstacles for global minima search over continuous spaces. For the semi-definite programming like (\ref{sdp_pro}), it has been known that the low-rank transformed problem (\ref{non-convex_pro}) has no local optima except the global ones if $p > \sqrt{2n}$ \cite{BM_smooth}. \par Thus, the main work is to discuss whether the proposed algorithm escapes strict saddle points and converges to global optimal solutions. At first, we provide the definition of unstable critical points. 
 Denote the update of $V$ generated by DSA as $V (t+1)= h_{SM}(V(t))$. 
\begin{definition}\label{cri_de2}
	Define unstable critical points as the set of critical points where the Jacobian of variable updating $h_{SM}(V)$ has at least a single eigenvalue with magnitude greater than one\upcite{saddle_escape},
	$$\mathcal{A}_g^*=\{V: h_{SM}(V)=V, \max_i\left|\lambda_i(Dh_{SM}(V))\right|>1\}.$$
\end{definition}
\par Following the work in \cite{saddle_escape}, we have the following property, which was investigated in \cite{wang2017mixing}.
\begin{lemma}\label{un_stable_fix}
	If $p>\sqrt{2n}$, each strict saddle point $V^*$ of the updating $h_{SM}$ is an unstable critical point, meaning $\mathcal{X}^*\subset \mathcal{A}_g^*$, where $\mathcal{X}^*$ is the set of strict saddle points.
\end{lemma}
\par\textbf{Proof:}
Consider the equivalent form in Lemma \ref{cmp_lem} of the proposed synchronous algorithm. It follows from the proof  in \cite{wang2017mixing} that  the Jacobi of proposed algorithm has eigenvalues containing those of the Jacobi of a standard Gauss-Seidel updating proposed in \cite{saddle_escape}. Based on the discussions of standard Gauss-Seidel updating in \cite{saddle_escape}, we obtain the desirable result.
$\hfill\blacksquare$
\par From Lemma \ref{un_stable_fix}, we deduce that all non-optimal critical pints are unstable fixed points. Next, we will prove that the updating $h_{SM}$ is a diffeomorphism, which is an invertible function that maps one differentiable manifold to another such that both the function and its inverse are smooth.
\begin{lemma}\label{diffeo_lemma}
	Under Assumptions \ref{SR_indice} and \ref{twoshare_assu}, the distributed synchronous updating $h_{SM}$ is a diffeomorphism.
\end{lemma}
\par\textbf{Proof:}
By the designed variable updating in Algorithm \ref{algo_sum} and Lemma 4.1, $h_{SM}$ is equivalent to
$$h_{SM}(V)=\left[\begin{matrix}\psi_n(\psi_{n-1}(\cdots\psi_1(V)))\end{matrix}\right]$$
where each column variable updating is defined as
$$ (\psi_i({V}))_{s=1\cdots n} =\left\{
\begin{array}{rcl}
&\frac{v_i-\Theta_{(i,i)} V M_{(:,i)} }{\left\|v_i-\Theta_{(i,i)}  V M_{(:,i)}\right\|} & \text{if} \quad s=i\\
&v_s' & \text{otherwise.}
\end{array} \right.$$
Because a composition of diffeomorphisms is still a diffeomorphism\upcite{inci2012regularity}, to prove this lemma, we only need to prove
 that $\psi_i(V)$ is a diffeomorphism for $i=1,\cdots,n$.
\par In \eqref{distri_mapping}, the step size takes a constant $\Theta_{(i,i)}\in(0,\frac{1}{\left\|M_{(i,:)}\right\|_1})$. Because $M$ is symmetric, it is equivalent to taking $\frac{1-\sigma}{\left\|M_{(:,i)}\right\|_1}$ for a constant $\sigma \in (0,1)$ and from the triangular inequality,
\begin{align*}
\left\| V M_{(:,i)}\right\|=\left\|\sum_{j=1}^n M_{(i,j)} v_j\right\|\leq \sum_{j=1}^n \left|M_{(i,j)}\right|\|v_j\|=\|M_{(i,:)}\|_1.
\end{align*}
Hence, $\left\|\Theta_{(i,i)} V M_{(:,i)}\right\|\leq 1-\sigma<1$. Thus, we have $\left\|v_i-\Theta_{(i,i)} V M_{(:,i)} \right\|\geq 1- \left\|\Theta_{(i,i)}  V M_{(:,i)}\right\|\geq \sigma>0$. Note that the function $\psi_i$ is only non-smooth at the point where the denominator term $\left\|v_i-\Theta_{(i,i)} V M_{(:,i)} \right\|=0$ and we have proved that the term is greater than $0$. Therefore, the function $\psi_i$ and its inverse function are valid and smooth. By the work in  \cite[Lemma C.2]{wang2017mixing}, $\psi_i$ is a diffeomorphism.  Since $h_{SM}(\cdot)$ is the composition of $\psi_i(\cdot)$s, the mapping of $h_{SM}(V)$ is also a diffeomorphism.
$\hfill\blacksquare$
\par Next, with the definition of global variable $V$ in \eqref{V_def}, we provide one equivalent form of objective function in optimization problem \eqref{non-convex_pro}, which will be used in the analysis of Lemma \ref{f_decrease}.
\begin{lemma}\label{ques_equal}
With the definition \eqref{V_def} and Assumptions \ref{SR_indice}, \ref{twoshare_assu}, the objective function of \eqref{non-convex_pro} at time $t$ $\sum_{\tilde{i}=1}^m \langle M^{\tilde{i}} ,V^{\tilde{i}'}(t)V^{\tilde{i}}(t)\rangle =\langle M,V(t)'V(t)\rangle$.
\end{lemma}
\par\textbf{Proof:}
The objective function of optimization problem \eqref{non-convex_pro} is $\sum_{\tilde{i}=1}^m \langle M^{\tilde{i}} ,V^{\tilde{i}'}V^{\tilde{i}}\rangle$. By the algorithm design, after one iteration $t$, $v_j^{\tilde{i}}(t)=v_j(t)$ for all $j\in J_{\tilde{i}}$, which holds because the coupling variables $v_j^{\tilde{i}}$ with $j\in S_{\tilde{i},{\rm par}(\tilde{i})}$ are equal to the uncoupling variables $v_j^{{\rm par}(\tilde{i})}$. Hence, the constraints \eqref{norm_con} and \eqref{share_con} in \eqref{non-convex_pro} hold. Then, we obtain
\begin{align*}
&\sum_{\tilde{i}=1}^m \langle M^{\tilde{i}} ,V^{\tilde{i}'}(t)V^{\tilde{i}}(t)\rangle\\
=&\sum_{\tilde{i}=1}^m \sum_{l,h\in J_{\tilde{i}}} M^{\tilde{i}}_{\{l,h\}} v^{\tilde{i}}_{l}(t)'v^{\tilde{i}}_{h}(t)\\
=&\sum_{\tilde{i}=1}^m \sum_{l,h\in J_{\tilde{i}}} M^{\tilde{i}}_{\{l,h\}} v_{l}(t)'v_{h}(t)\\
=&\sum_{l,h\in J_{\tilde{i}}} M_{(l,h)} v_{l}(t)'v_{h}(t)\\
=&\langle M,V(t)'V(t)\rangle,
\end{align*} 
where the second to last equation holds because $M^{\tilde{i}}_{\{l,h\}}=\overline{M}^{\tilde{i}}_{(l,h)}$ and $\sum_{\tilde{i}=1}^m \overline{M}^{\tilde{i}}=M$.
$\hfill\blacksquare$
\par With Lemma \ref{ques_equal}, before updating $v_i$, all variable $v_j$ except for $v_i$ are given and fixed, then the global function is rewritten as
\begin{align}\label{trans_pro}
f(V)&=\left<M, V'V\right>\notag\\&=\sum_{i=1}^n \sum_{j=1}^n M_{(i,j)} v_i'v_j \notag\\
&=2v_i'\mathbf{g}_i+{\rm constant},
\end{align}
where the last equation holds since the matrix $M$ is symmetric. Note that $ \mathbf{g}_i$ is defined in Lemma \ref{cmp_lem} and is independent of $v_i$ because $M_{(i,i)}=0$. Then, we have the following Lemma stating the monotonous decreasing property of the global function value generated by the proposed synchronous algorithm.
\begin{lemma}\label{f_decrease}
	For the proposed synchronous algorithm with step size $\Theta$, let ${V}(t+1)=h_{SM}(V(t))$. Under Assumptions \ref{SR_indice} and \ref{twoshare_assu}, we have
	\begin{align}\label{f_f}
	f(V(t))-f(V(t+1))=\sum_{i=1}^n \frac{1+y_i(t)}{\Theta_{(i,i)}} \left\|v_i(t)-v_i(t+1)\right\|^2,
	\end{align}
	where $y_i(t)=\left\|v_i(t)-\Theta_{(i,i)}\mathbf{g}_i(t)\right\|$ and  $\mathbf{g}_i(t)=\sum_{l<i} M_{(i,l)} v_l(t+1)+\sum_{l>i} M_{(i,l)} v_l(t)$.
\end{lemma}
\par\textbf{Proof:}
By (\ref{trans_pro}), the function difference after updating $v_i(t)$ to $v_i(t+1)$ is $2 \mathbf{g}_i'(v_i(t)-{v}_i(t+1))$. Then, by the updating in (\ref{distri_mapping}), ${v}_i(t+1)=(v_i(t)-\Theta_{(i,i)}\mathbf{g}_i(t))/ y_i(t)$, we have
\begin{align}\label{f_de}
&2 \mathbf{g}_i(t)'(v_i(t)-{v}_i(t+1))\notag\\
=& 2(\mathbf{g}_i(t)+\frac{v_i(t)-\Theta_{(i,i)}\mathbf{g}_i(t)}{\Theta_{(i,i)}})'(v_i(t)-{v}_i(t+1))\notag\\
&-2(\frac{v_i(t)-\Theta_{(i,i)}\mathbf{g}_i(t)}{\Theta_{(i,i)}})'(v_i(t)-{v}_i(t+1))\notag\\
=&2 \frac{1}{\Theta_{(i,i)}}v_i(t)'(v_i(t)-{v}_i(t+1))\notag\\
&-2\frac{y_i(t)}{\Theta_{(i,i)}}{v}_i(t+1)'(v_i(t)-{v}_i(t+1))\notag\\
=&\frac{1+y_i(t)}{\Theta_{(i,i)}} 2 (1-v_i(t)'{v}_i(t+1))\notag\\
=&\frac{1+y_i(t)}{\Theta_{(i,i)}} \left\|v_i(t)-{v}_i(t+1)\right\|^2,
\end{align}
where the third equality holds due to the condition $\|v_i\|=1$.
Then, the result holds from summing the above equation over $i=1,\cdots,n$.
$\hfill\blacksquare$
\par Now, we are ready to prove the result in Theorem \ref{syn_theo}.
\par \textbf{Proof of Theorem \ref{syn_theo}:}
Assume Assumptions \ref{SR_indice} and \ref{twoshare_assu} hold. From Lemma $\ref{un_stable_fix}$ and the nonexistence of local optima, all non-optimal critical pints are unstable fixed points. Recall that  $h_{SM}$ is a diffeomorphism by Lemma \ref{diffeo_lemma} and non-optimal critical points of $h_{SM}(\cdot)$ are unstable fixed points by Lemma \ref{un_stable_fix}. It follows from the center-stable manifold theorem (Theorem III.5 of \cite{shub_book}) that the proposed algorithm escapes all non-optimal critical points almost surely under random initialization. By Lemma \ref{f_decrease} and the fact that  $\frac{1+y_i(t)}{\Theta_{(i,i)}}$ in (\ref{f_f}) is always positive over iterations, the objective function value is strictly decreasing. Because the objective function value generated by the proposed algorithm is strictly decreasing and the objective value is lower bounded, the generated variables converge to the set of  first-order critical points. {Thus, the almost sure divergence from the non-optimal critical points and the convergence to critical points imply that $v_j^{\tilde{i}}, j\in J_{\tilde{i}}$ in the updating $h_{SM}$ converges to corresponding column of global optimal solutions of (\ref{non-convex_pro}) almost surely under random initialization.
\par Next, we show that the result of this theorem holds if Assumptions \ref{SR_indice} and \ref{twoshare_assu} are removed. If Assumption \ref{SR_indice} does not hold, the indices of variables can be rearranged manually such that the indices in uncoupling set $\mathcal{R}_{\tilde{i},{\rm par}(\tilde{i})}$ are smaller than the indices in coupling set $\mathcal{S}_{\tilde{i},{\rm par}(\tilde{i})}$ for each agent $\tilde{i}$. Hence, the above analysis still holds without Assumptions \ref{SR_indice}. If Assumption \ref{twoshare_assu} does not hold, in \eqref{Slgeqj} of Lemma \ref{cmp_lem}, for each variable $v_l^{\tilde{i}}, l\in \mathcal{S}_{{\rm ch}(\tilde{i}), \tilde{i}}\cap \mathcal{S}_{\tilde{i},{\rm par}(\tilde{i})}$, there must be a parent agent $\tilde{j}$ such that $l\in \mathcal{R}_{\tilde{j},{\rm par}(\tilde{j})}$ by algorithm design. Then, $v_l^{\tilde{i}}(t)=v_l^{\tilde{j}}(t)=v_l(t)$. Thus, the analysis in \eqref{Slgeqj} analogously holds and the rest of theoretical deductive is true.}
$\hfill\blacksquare$
\subsection{Convergence analysis for asynchronous algorithm }\label{ASM_sec}
\par {For distributed asynchronous algorithm \ref{algo_asyn}}, let $T^{\tilde{i}}$ be the set of times at which agent $\tilde{i}$ updates variable $V^{\tilde{i}}$. In addition, agent $\tilde{i}$ may not have access to the most recent value of other agents' variables. To collect communicated information from neigbors, define a set $\mathcal{J}_{\tilde{i}}$ as the union $J_{\tilde{i}}\cup \mathcal{R}_{\tilde{r}\in {\rm ch}(\tilde{i}),\tilde{i}}$. Thus, define one possibly outdated variable of agent $\tilde{i}$, {$\mathbf{V}^{\tilde{i}}(t)\in \mathbb{R}^{p\times |\mathcal{J}_{\tilde{i}}|}$}, as
\begin{align}\label{Vi_def}
\mathbf{V}^{\tilde{i}}(t)=\big[v_{s_1}^{\tilde{i}}(\tau_{s_1}^{\tilde{i}}(t)),\cdots,v_{s_{|\mathcal{J}_{\tilde{i}}|}}^{\tilde{i}}(\tau_{s_{|\mathcal{J}_{\tilde{i}}|}}^{\tilde{i}}(t))\big],
\end{align}
where $\{s_1,\cdots,s_{|\mathcal{J}_{\tilde{i}}|}\}=\mathcal{J}_{\tilde{i}}$, $\tau_{{j}}^{\tilde{i}}(t)$ is assumed to satisfy the condition (\ref{taulim}). Recall that, if $s_j\in \mathcal{R}_{\tilde{i},{\rm par}(\tilde{i})}$, then $\tau_{s_j}^{\tilde{i}}(t)=t$. $\mathbf{V}^{\tilde{i}}\in \mathbb{R}^{p\times |\mathcal{J}_{\tilde{i}}|}$ collects time-delayed transmitted information from children and parents. 

\par With additional assumption in Assumption \ref{asyn_assump} that each element of global coefficient matrix $M$ is only accessible to one agent, the relationship $M^{\tilde{i}}_{\{l,k\}}=M_{(l,k)}$ holds for all agent $\tilde{i}$. Then, the updating proposed in DAA is rewritten as a compact form as shown in the following lemma, where each column variable $v_j$ of global variabel $V$ defined in \eqref{V_def} is expressed by outdated transmitted information.
\begin{lemma}\label{asyn_comp}
	Under Assumption \ref{asyn_assump}, each column variable $v_j$ generated by DAA with index $j\in \mathcal{R}_{\tilde{i},{\rm par}(\tilde{i})}$ is equivalent to
	\begin{equation}
	\left\{
	\begin{aligned}
	&v_j(t+1)={\rm normal}(v_j(t)-\Theta_{(j,j)}\mathbf{h}_j(\mathbf{V}^{\tilde{i}}(t)), \quad t\in T^{\tilde{i}},\\
	&v_j(t+1)=v_j(t),\quad t\notin T^{\tilde{i}}
	\end{aligned}
	\right.
	\end{equation}
	where $\mathbf{h}_j(\mathbf{V}^{\tilde{i}}(t))=  \big[\sum_{l\in \mathcal{R}_{\tilde{i},{\rm par}(\tilde{i})}} M_{(j,l)} {v}_l(t)+\sum_{l\in \mathcal{S}_{\tilde{i},{\rm par}(\tilde{i})}} M_{(j,l)} {v}_l(\tau_{l}^{\tilde{i}}(t))+\sum_{\tilde{r}\in {\rm ch}(\tilde{i})}\sum_{l\in \mathcal{R}_{\tilde{r},\tilde{i}}} M_{(j,l)} v_l(\tau_{l}^{\tilde{i}}(t))\big]$, $\Theta_{(j,j)}=\theta_j^{\tilde{i}}$ was defined in (\ref{asyn_step}).
\end{lemma}
\par\textbf{Proof:}
By (\ref{v_updateasm}) in DAA, agent $\tilde{i}$ updates local uncoupling variable $v_j$ with index in $j\in \mathcal{R}_{\tilde{i},{\rm par}(\tilde{i})}$
according to
\begin{align}\label{vrasyn}
&v_j^{\tilde{i}}(t+1)\notag\\=&{\rm normal}\Big(v_j^{\tilde{i}}(t)-\theta_j^{\tilde{i}} \big[p_j^{\tilde{i}}+\sum_{\tilde{r}\in {\rm ch}(\tilde{i})}\varpi_j^{\tilde{r},\tilde{i}}(\tau^{\tilde{i}}(t))\big]\Big)\notag\\
=&{\rm normal}\Big(v_j^{\tilde{i}}(t)-\theta_j^{\tilde{i}} \big[\sum_{l\in J_{\tilde{i}}} M_{\{j,l\}}^{\tilde{i}} {v}_l^{\tilde{i}}(t)\notag\\
&+\sum_{\tilde{r}\in {\rm ch}(\tilde{i})}\sum_{l\in J_{\tilde{r}}} M_{\{j,l\}}^{\tilde{r}} v_l^{\tilde{r}}(\tau_{l}^{\tilde{i}}(t))\big]\Big)\notag\\
=&{\rm normal}\Big(v_j^{\tilde{i}}(t)-\theta_j^{\tilde{i}} \big[\sum_{l\in \mathcal{R}_{\tilde{i},{\rm par}(\tilde{i})}} M_{(j,l)} {v}_l^{\tilde{i}}(t)\notag\\
&+\sum_{l\in \mathcal{S}_{\tilde{i},{\rm par}(\tilde{i})}} M_{(j,l)} {v}_l^{{\rm par}({\tilde{i}})}(\tau_{l}^{\tilde{i}}(t))\notag\\
&+\sum_{\tilde{r}\in {\rm ch}(\tilde{i})}\sum_{l\in \mathcal{R}_{\tilde{r},\tilde{i}}} M_{(j,l)} v_l^{\tilde{r}}(\tau_{l}^{\tilde{i}}(t))\big]\Big)
\end{align}
where the last equality holds because $M_{\{j,l\}}^{\tilde{r}}=0$, for $j\in \mathcal{R}_{\tilde{i},{\rm par}(\tilde{i})}$, $l\in \mathcal{S}_{\tilde{r},\tilde{i}}$, by the assumption that each element of $M$ is only accessible to one agent. 
\par For the second term $\sum_{l\in \mathcal{S}_{\tilde{i},{\rm par}(\tilde{i})}} M_{(j,l)} {v}_l^{{\rm par}({\tilde{i}})}(\tau_{l}^{\tilde{i}}(t))$ of (\ref{vrasyn}), similarly to the discussions of (\ref{Slgeqj}) in the synchronous case, we obtain
\begin{align}\label{secondterm}
\!\sum_{l\in \mathcal{S}_{\tilde{i},{\rm par}(\tilde{i})}}\! M_{(j,l)} {v}_l^{{\rm par}({\tilde{i}})}\!(\tau_{l}^{\tilde{i}}(t))\!=\!\sum_{l\in \mathcal{S}_{\tilde{i},{\rm par}(\tilde{i})}} \! M_{(j,l)} {v}_l(\tau_{l}^{\tilde{i}}(t)).
\end{align}

\par Then, substituting (\ref{secondterm}) to (\ref{vrasyn}), we have
\begin{align}\label{ljasy}
&v_j(t+1)\notag\\
=&{\rm normal}\Big(v_j(t)-\Theta_{(j,j)} \big[\sum_{l\in \mathcal{R}_{\tilde{i},{\rm par}(\tilde{i})}} M_{(j,l)} {v}_l^{\tilde{i}}(t)\notag\\
&+\sum_{l\in \mathcal{S}_{\tilde{i},{\rm par}(\tilde{i})}} M_{(j,l)} {v}_l(\tau_{l}^{\tilde{i}}(t))\notag\\
&+\sum_{\tilde{r}\in {\rm ch}(\tilde{i})}\sum_{l\in \mathcal{R}_{\tilde{r},\tilde{i}}} M_{(j,l)} v_l^{\tilde{r}}(\tau_{l}^{\tilde{i}}(t))\big]\Big)\notag\\
=&{\rm normal}\Big(v_j(t)-\Theta_{(j,j)}\big[\sum_{l\in \mathcal{R}_{\tilde{i},{\rm par}(\tilde{i})}} M_{(j,l)} {v}_l(t)\notag\\
&+\sum_{l\in \mathcal{S}_{\tilde{i},{\rm par}(\tilde{i})}} M_{(j,l)} {v}_l(\tau_{l}^{\tilde{i}}(t))\notag\\
&+\sum_{\tilde{r}\in {\rm ch}(\tilde{i})}\sum_{l\in \mathcal{R}_{\tilde{r},\tilde{i}}} M_{(j,l)} v_l(\tau_{l}^{\tilde{i}}(t))\big]\Big),\notag\\
&={\rm normal}(v_j(t)-\Theta_{(j,j)} \mathbf{h}_j(\mathbf{V}^{\tilde{i}}(t))), \quad t\in T^{\tilde{i}},
\end{align}
where {$\mathbf{h}_j(\mathbf{V}^{\tilde{i}}(t))= \big[\sum_{l\in \mathcal{R}_{\tilde{i},{\rm par}(\tilde{i})}} M_{(j,l)} {v}_l(t)+\sum_{l\in \mathcal{S}_{\tilde{i},{\rm par}(\tilde{i})}} M_{(j,l)} {v}_l(\tau_{l}^{\tilde{i}}(t))+\sum_{\tilde{r}\in {\rm ch}(\tilde{i})}\sum_{l\in \mathcal{R}_{\tilde{r},\tilde{i}}} M_{(j,l)} v_l(\tau_{l}^{\tilde{i}}(t))\big]$}, $\Theta_{(j,j)}=\theta_j^{\tilde{i}}$ for $j\in \mathcal{R}_{\tilde{i},{\rm par}(\tilde{i})}$. For $t\notin T^{\tilde{i}}$, $v_j(t+1)=v_j(t)$ holds naturally.
$\hfill\blacksquare$
\par Define the updating direction $s_j$ ($j\in \mathcal{R}_{\tilde{i},{\rm par}(\tilde{i})}$) as, 
\begin{align}\label{s_de}
s_j(t)&=\frac{1}{\Theta_{(j,j)}}(v_j(t+1)-v_j(t)).
\end{align}
If $t\in T^{\tilde{i}}$, 
\begin{align}
s_j(t)=\frac{1}{\Theta_{(j,j)}}\!\Big(\!{\rm normal}\big(v_j(t)-\Theta_{(j,j)}\mathbf{h}_j(\mathbf{V}^{\tilde{i}}(t))\big)-v_j(t)\!\Big),
\end{align}
 and if $t\notin T^{\tilde{i}}$, $s_j(t)=0$.
\par In the following lemma, we present a vital descent property of local variable iteration, which will be used in the proof of Theorem \ref{asyn_theo}.
\begin{lemma}\label{coor_descent}
	Suppose Assumption \ref{asyn_assump} holds. For any agent $\tilde{i}$ and time $t$, we have
	\begin{align}\label{sineq}
	s_j(t)'\mathbf{h}_j(\mathbf{V}^{\tilde{i}}(t))\leq -\|s_j(t)\|^2, \quad j\in \mathcal{R}_{\tilde{i},{\rm par}(\tilde{i})}.
	\end{align}
\end{lemma}
\textbf{Proof:}
If $t\notin T^{\tilde{i}}$, the inequality (\ref{sineq}) is true since both sides are zero. If $t\in T^{\tilde{i}}$, by the definition of $s_j(t)$ in (\ref{s_de}) and by Lemma \ref{asyn_comp}, $v_j(t+1)=(v_j(t)-\Theta_{(j,j)}\mathbf{h}_j(\mathbf{V}^{\tilde{i}}(t)))/y_j(t)$ for $j\in \mathcal{R}_{\tilde{i},{\rm par}(\tilde{i})}$, where $y_j(t)=\|v_j(t)-\Theta_{(j,j)}\mathbf{h}_j(\mathbf{V}^{\tilde{i}}(t))\|$, we have
\begin{align*}
&s_j(t)'\mathbf{h}_j(\mathbf{V}^{\tilde{i}}(t))\\
=& \frac{-1}{\Theta_{(j,j)}}(v_j(t)-v_j(t+1))'\mathbf{h}_j(\mathbf{V}^{\tilde{i}}(t))\\
=& \frac{-1}{\Theta_{(j,j)}} \bigg[\Big(\mathbf{h}_j(\mathbf{V}^{\tilde{i}}(t))\\
&+\frac{v_j(t)-\Theta_{(j,j)}\mathbf{h}_j(\mathbf{V}^{\tilde{i}}(t))}{\Theta_{(j,j)}}\Big)'(v_j(t)-v_j(t+1))\\
&-\Big(\frac{v_j(t)-\Theta_{(j,j)}\mathbf{h}_j(\mathbf{V}^{\tilde{i}}(t))}{\Theta_{(j,j)}}\Big)' (v_j(t)-v_j(t+1)) \bigg]\\
=&  \frac{-1}{\Theta_{(j,j)}} \bigg[\frac{1}{\Theta_{(j,j)}} v_j(t)'(v_j(t)-v_j(t+1))\\
&-\frac{y_j(t)}{\Theta_{(j,j)}} v_j(t+1)'(v_j(t)-v_j(t+1))\bigg]\\
=& -\frac{1+y_j(t)}{\Theta_{(j,j)}^2}\|v_j(t+1)-v_j(t)\|^2\\
=&-(1+y_j(t))\|s_j(t)\|^2\\
\leq &-\|s_j(t)\|^2,
\end{align*}
where the last inequality holds because $y_j$ is non-negative.
$\hfill\blacksquare$
\par Making use of Lemma \ref{coor_descent}, we discuss the relationship of global variable $V$ and local variables $\mathbf{V}^{\tilde{i}}$, and the gradient of objective function at the point $V(t_k)$ when $k \to \infty$ in Theorem \ref{asyn_theo}. Before discussions, it should be noted that $\mathbf{h}_j(\mathbf{V}^{\tilde{i}}(t))$ defined in Lemma \ref{asyn_comp} is exactly the gradient of global function $f$ with respect to column variable $v_j(t)$, where $j\in \mathcal{R}_{\tilde{i},{\rm par}(\tilde{i})}$, so for convenience, we use $\mathbf{h}_j(\mathbf{V}^{\tilde{i}}(t))$ in the following analysis. The proof follows the studies of gradient-like optimization algorithms in Proposition 5.1, section 7, \cite{paral_distri_book}.
\par \textbf{Proof of Theorem \ref{asyn_theo}:} We follow the proof of Proposition 5.1 in  \cite{paral_distri_book}. By Assumption \ref{pro_assump} (2), the objective function of \eqref{non-convex_pro} satisfies $\sum_{\tilde{i}=1}^m f_{\tilde{i}}\geq 0$. In addition, with the analysis in Lemma \ref{coor_descent}, we have, for $j\in \{1,\cdots,n\}$,
\begin{align*}
s_j(t)'\mathbf{h}_j(\mathbf{V}^{\tilde{i}}(t))&\geq -\|s_j(t)\|\|\mathbf{h}_j(\mathbf{V}^{\tilde{i}}(t))\|\\
-(1+y_j(t))\|s_j(t)\|^2 &\geq -\|s_j(t)\|\|\mathbf{h}_j(\mathbf{V}^{\tilde{i}}(t))\|\\
(1+y_j(t))\|s_j(t)\|^2 &\leq \|s_j(t)\|\|\mathbf{h}_j(\mathbf{V}^{\tilde{i}}(t))\|\\
\|s_j(t)\|&\leq \frac{1}{1+y_j(t)}\|\mathbf{h}_j(\mathbf{V}^{\tilde{i}}(t))\|,
\end{align*}
where the last inequality holds because $y_j$ is non-negative. Then, there is a positive constant $K_3=1$ such that $\|s_j(t)\|\leq K_3 \|\mathbf{h}_j(\mathbf{V}^{\tilde{i}}(t))\|$. 
What's more, {with the block-descent property in Lemma \ref{coor_descent},} the assumptions in Proposition 5.1 \cite{paral_distri_book} hold, where the  product is replaced by inner product of vectors and the term $\left|s_j(t)\right|$ is replaced by $\|s_j(t)\|$. Then, by a similar analysis as Proposition 5.1 in \cite{paral_distri_book}, we obtain that 
\begin{align}\label{s_lim}
\lim_{t \to \infty} s_j(t)=0,
\end{align}
for each $j\in\{1,\cdots,n\}$. In addition, by (\ref{s_de}), we obtain
\begin{align}
&\lim_{t\to \infty}\|V(t+1)-V(t)\|=0.\label{Vgap}
\end{align}
Then, consider the boundeness of $\|v_j^{\tilde{i}}(t)-v_j(t)\|$.
\begin{align}\label{vv}
\|v_j^{\tilde{i}}(t)-v_j(t)\|=&\|v_j^{\tilde{i}}(\tau_j^{\tilde{i}}(t))-v_j(t)\|\notag\\
=&\Theta_{(j,j)}\|\sum_{\tau=\tau_{j}^{\tilde{i}}(t)}^{t-1}s_j(\tau)\|\notag\\
\leq &\Theta_{(j,j)}\sum_{\tau=t-B}^{t-1}\|s_j(\tau)\|.
\end{align}
With \eqref{s_lim} and \eqref{vv}, we also obtain
\begin{align}
&\lim_{t\to \infty}\|v_j^{\tilde{i}}(t)-v_j(t)\|=0, \ \forall j\in \mathcal{R}_{\tilde{i},{\rm par}(\tilde{i})}.\label{ViVgap}
\end{align}
In addition, for each $j\in \{1,\cdots,n\}$, there is a unique agent $\tilde{i}$ such that $j\in \mathcal{R}_{\tilde{i},{\rm par}(\tilde{i})}$ and the equation \eqref{ViVgap} holds.
\par Define $V^*$ as a limit point of $V(t)$ and $\{t_k\}$ as a sequence such that $\lim_{k\to\infty}V(t_k)=V^*$. Let $\tau_k$ be such that $\left|t_k-\tau_k\right|\leq B$ and $\tau_k\in T^{\tilde{i}}$. Then, by equations (\ref{Vgap}) and (\ref{ViVgap}), $v_j(\tau_k)$  converges to $v_j^*$, which is $j$th column of $V^*$, and $v_j^{\tilde{i}}(\tau_k)$ converges to $v_j^*$ such that $\mathbf{V}^{\tilde{i}}(\tau_k)$ converges. Then, we have
\begin{align*}
&\lim_{k\to \infty}({\rm normal}(v_j^{\tilde{i}}(\tau_k)-\Theta_{(j,j)}\mathbf{h}_j(\mathbf{V}^{\tilde{i}}(\tau_k))-v_j^{\tilde{i}}(\tau_k))\\
=&\lim_{k\to \infty}\Theta_{(j,j)}s_j(\tau_k)=0.
\end{align*}
Since it holds for each $j\in\{1,\cdots,n\}$, we get the desired result.
$\hfill\blacksquare$

\section{Simulation}\label{simulation}
In this section, numerical tests and large-scale image segmentation application are presented to show the efficiency of the proposed distributed algorithms.
\par
\begin{figure}
	\centering
	\includegraphics[width=6cm]{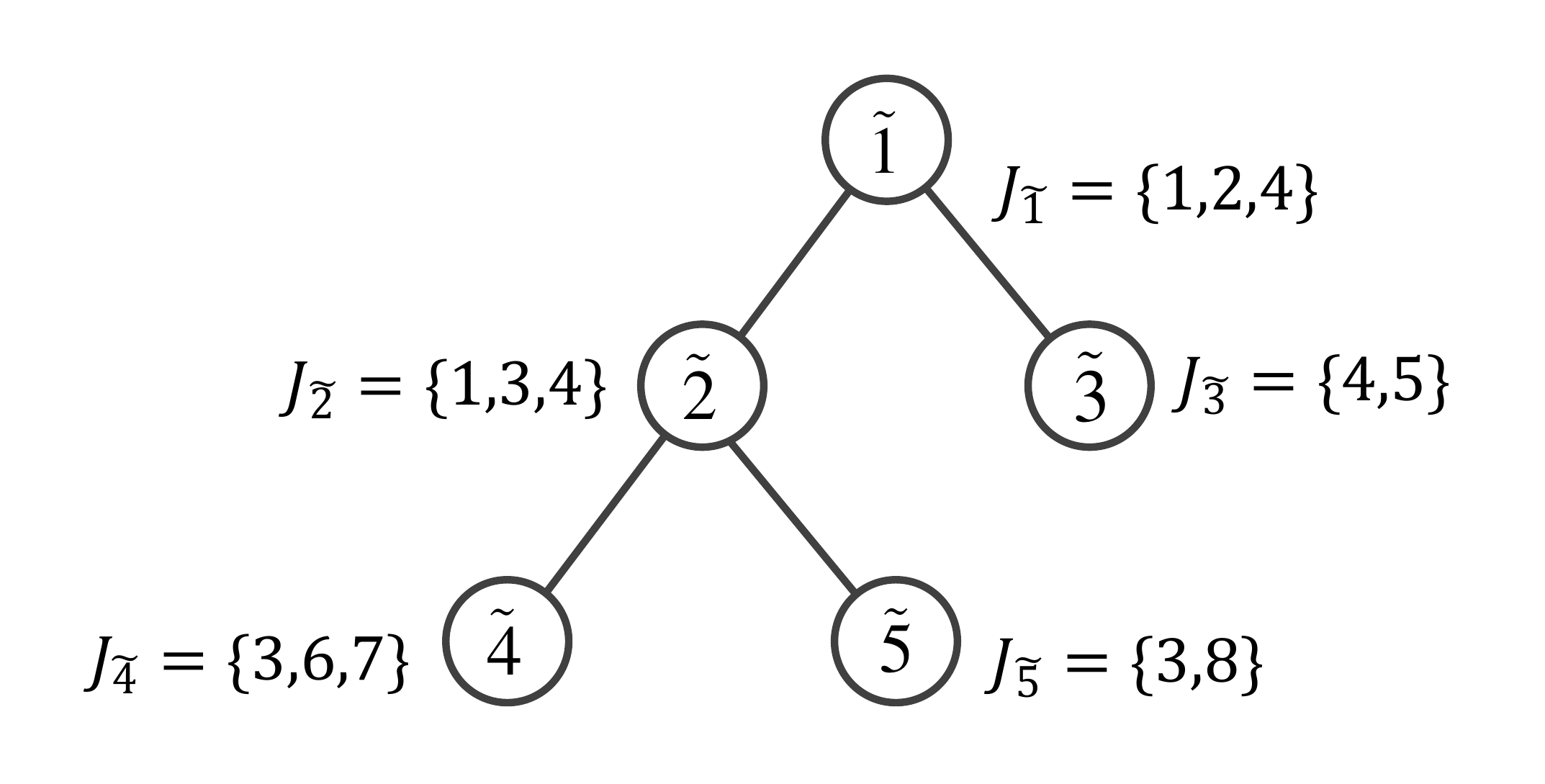}\\
	\caption{Clique tree of simulation problem}\label{simu_graph}
	\centering
\end{figure}

 \textit{Example 1:} We present one special sparse coupling numerical optimization problem, which has been investigated in \cite{dis_pd,distri_sed_tac,fast_dis} and of which the corresponding connected graph is one clique tree. More detailed information of clique trees can be found in \cite{dis_pd}. Here we only introduce some brief concepts and focus on the discussions about the numerical convergence performance of proposed algorithms. For the coupling optimization (\ref{sdp_pro}), we assume that the dimension of global matrix variable $X$ is $n=8$, the number of local functions is $N=6$, the corresponding dependent element indices set are $C_1=\{1,3\}$, $C_2=\{1,2,4\}$, $C_3=\{4,5\}$, $C_4=\{3,4\}$, $C_5=\{3,6,7\}$, $C_6=\{3,8\}$. By the clique tree transformations in \cite{distri_sed_tac}, the corresponding clique tree owns five agents, shown in Fig.\ref{simu_graph}. It shows that the number of agents in corresponding problem (\ref{non-convex_pro}) is $m=5$, which implies that one agent has multiple local functions. The local functions assigned to $i$th agent are denoted by a function set $\phi_{\tilde{i}}$. Then, the function sets of the clique tree are $\phi_{\tilde{1}}=\{f_2\}$, $\phi_{\tilde{2}}=\{f_1,f_4\}$, $\phi_{\tilde{3}}=\{f_3\}$, $\phi_{\tilde{4}}=\{f_5\}$, $\phi_{\tilde{5}}=\{f_6\}$. In addition, the ordered index sets are $J_{\tilde{1}}=\{1,2,4\}$, $J_{\tilde{2}}=\{1,3,4\}$, $J_{\tilde{3}}=\{4,5\}$, $J_{\tilde{4}}=\{3,6,7\}$, $J_{\tilde{5}}=\{3,8\}$. More specifically, we provide the decomposed diagram of sparse coefficient matrix $M$ over five different agents as following. Elements with different colors are assigned to different agents.
$$\left[
\begin{array}{c|c|c|c|c|c|c|c}
W_{11} &{\textcolor{red}{W_{12}}}&{\textcolor{blue}{W_{13}}}&{\textcolor{red}{W_{14}}}&0&0&0&0\\ \cdashline{1-8}[0.8pt/2pt]
{\textcolor{red}{W_{21}}} &W_{22}&0&{\textcolor{red}{W_{24}}}&0&0&0&0\\ \cdashline{1-8}[0.8pt/2pt]
{\textcolor{blue}{W_{31}}} &0&W_{33}&{\textcolor{blue}{W_{34}}}&0&{\textcolor{cyan}{W_{36}}}&{\textcolor{cyan}{W_{37}}}&{\textcolor{magenta}{W_{38}}}\\ \cdashline{1-8}[0.8pt/2pt]
{\textcolor{red}{W_{41}}} &{\textcolor{red}{W_{42}}}&{\textcolor{blue}{W_{43}}}&W_{44}&{\textcolor{green}{W_{45}}}&0&0&0\\ \cdashline{1-8}[0.8pt/2pt]
0 &0&0&{\textcolor{green}{W_{54}}}&W_{55}&0&0&0\\ \cdashline{1-8}[0.8pt/2pt]
0 &0&{\textcolor{cyan}{W_{63}}}&0&0&0&{\textcolor{cyan}{W_{67}}}&0\\ \cdashline{1-8}[0.8pt/2pt]
0 &0&{\textcolor{cyan}{W_{73}}}&0&0&{\textcolor{cyan}{W_{76}}}&W_{77}&0\\ \cdashline{1-8}[0.8pt/2pt]
0 &0&{\textcolor{magenta}{W_{83}}}&0&0&0&0&W_{88}
\end{array}
\right].$$
\par Since the diagonal elements will not influence the optimization result, they are assigned to any agent without effect. Let $f^*$ be the optimal function value of the optimization problem (\ref{sdp_pro}), which is solved by the solver YALMIP when the dimension $n$ is not too large. 
\par
\begin{figure}
	\centering
	\subfigure[The trajectories of $f$ by DSA, DAA and SDPLR]{
		\includegraphics[width=8cm]{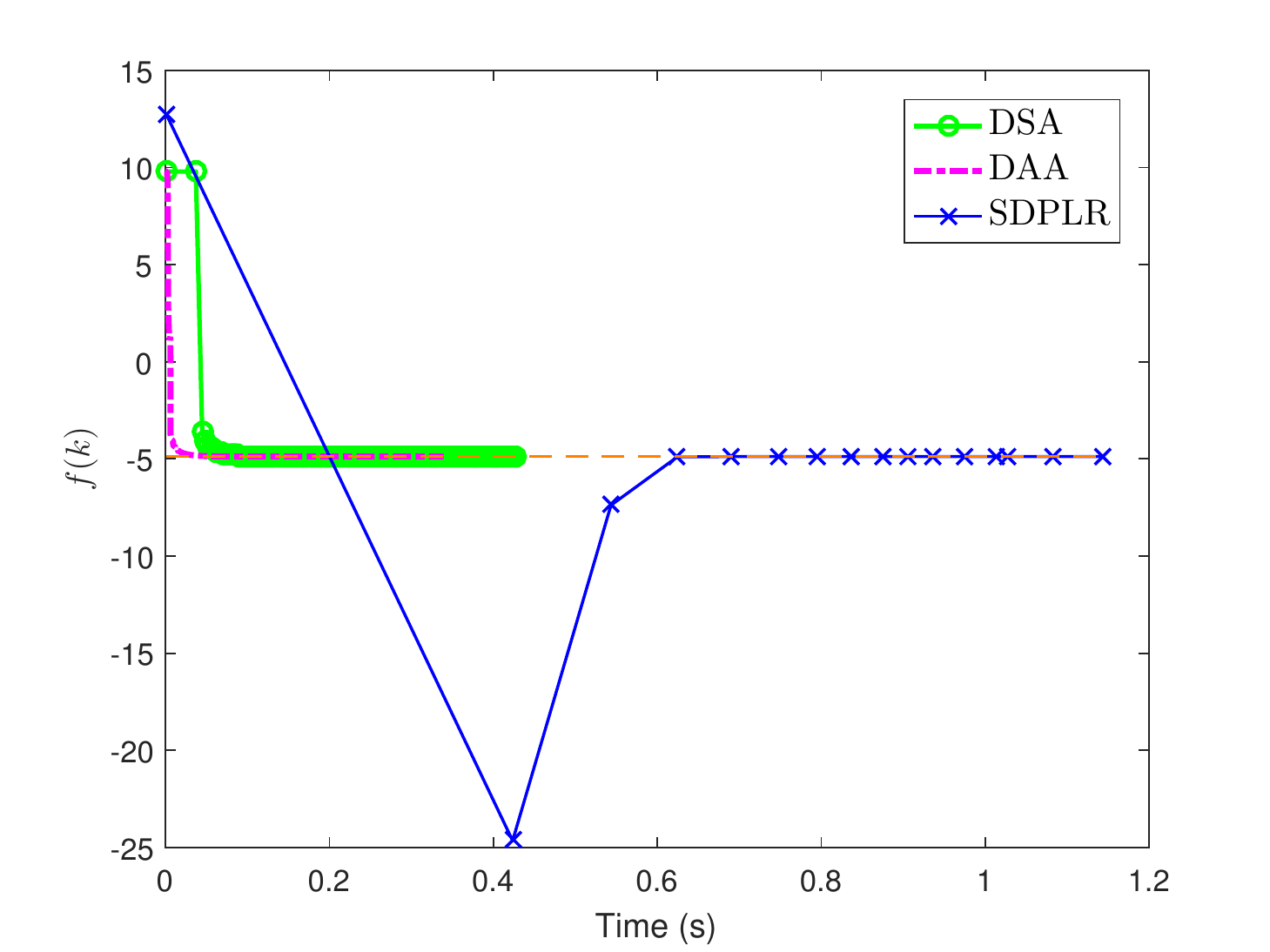}
		\label{compare_fig}
	}
	
	\subfigure[The trajectories of $\left\|X(k+1)-X(k)\right\|_F$ along time]{
		\includegraphics[width=8cm]{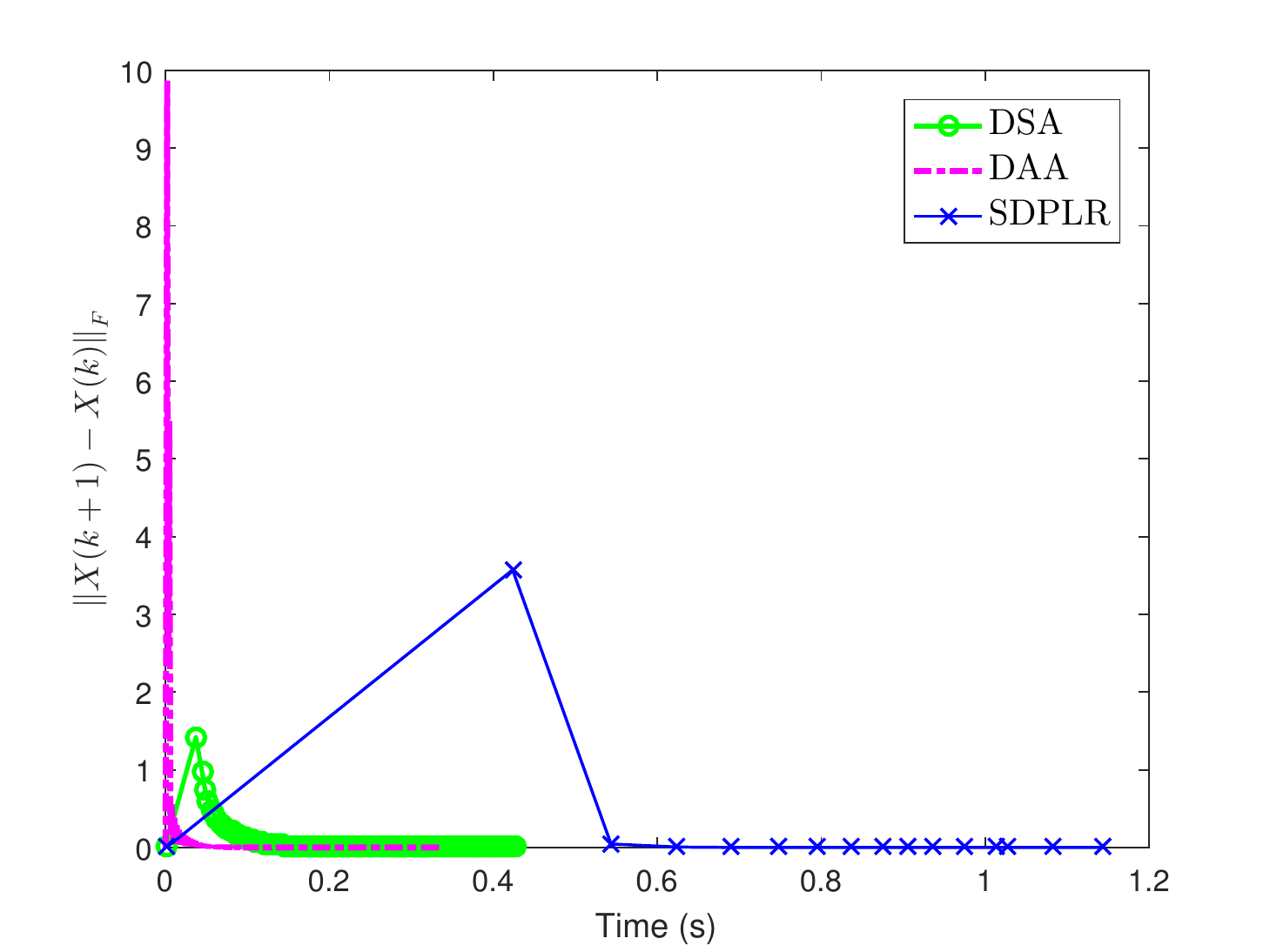}
		\label{gapx}
	}
	\subfigure[The trajectories of gradient of $f$ along time]{
		\includegraphics[width=8cm]{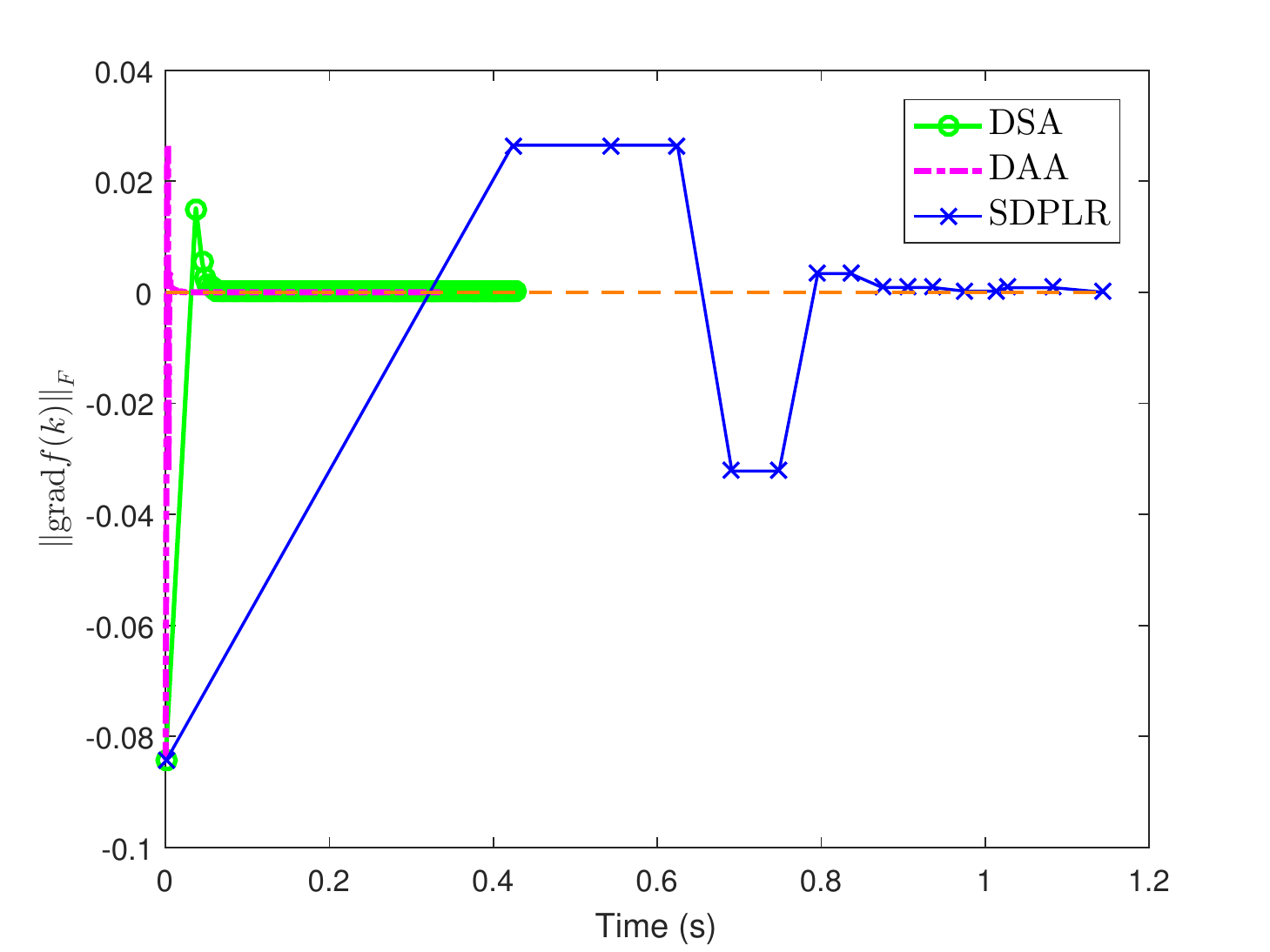}
		\label{gradf}
	} 	
	\caption{Convergent trajectories generated by proposed algorithms and SDPLR}
\end{figure}

 $1)$ We use the proposed distributed synchronous, asynchronous algorithms and the centralized algorithm SDPLR \cite{BM_SPLR}, which are all coded by MATLAB, to solve the sparse optimization problem. The simulation results are shown in Figs. \ref{compare_fig}-\ref{gradf}. The original global variable value is $X(k)=V(k)'V(k)$. In Fig. \ref{compare_fig}, the proposed algorithms and SDPLR all converge to the optimal function value $f^*$, which is calculated by the solver YALMIP. It is observed that distributed algorithms converge much faster than the SDPLR algorithm for the MAXCUT problem. In Fig. \ref{gapx}, the trajectories of $\left\|X(k+1)-X(k)\right\|_F$, where the trajectories $\{X(k)\}$ are generated by the proposed algorithms DSA and DAA respectively, are shown to converge to zeros. It shows that varaible $X(k)$ converges to one limiting point.
The norm of the Riemannian gradient of $f$ is defined as $\left\|{\rm grad}f\right\|_F^2=\sum_{i=1}^n(\left\|{g}_i\right\|^2-\left<v_i,{g}_i\right>^2)$, where $\|\cdot\|_F$ represents the Frobenius norm of a matrix, and $g_i=\sum_{j=1}^n M_{(i,j)}v_j$. The trajectories of $\left\|{\rm grad} f(k)\right\|_F$ are shown in Fig. \ref{gradf} and converge to zeros, which implies that the generated sequences $\{V(k)\}$ converge to critical points. In addition, by the numerical experiments, for most cases, we have observed that the asynchronous algorithm converges faster than the synchronous algorithm. In addition, in the next simulation, we provide quantitative comparisons of convergence rates between DAA and DSA.
\par $2)$ In order to compare the performance of the proposed distributed algorithms with the inspired centralized algorithm, which is proposed in \cite{wang2017mixing}, we make use of MPICH distributed model, which is a high-performance message passing interface, to develop a multi-processers environment on one computer with a Core(TM) I5-8250U CPU, 1.6GHz. Both centralized algorithm and distributed algorithms are coded by C language. For the distributed algorithms, we use five processes to deal with the optimization problem.
\par We provide two experiments with different dimensions and collect the number of iterations and executive time of different algorithms. In each experiment, the stop criterion of iterations reaches an expected error between the function value $f(k)$ and optimal value $f^*$. In addition, we use $sn$ to denote the number of shared variables over the multi-agent network, e.g., the $sn$ of network shown in Fig. \ref{simu_graph} is $5$. The executed time comparisons of centralized algorithm and distributed algorithms are listed in following table \ref{com_table}.
\par By the comparative test, the distributed synchronous and asynchronous algorithms both converge faster than the centralized Mixing algorithm. As the dimension of problem increases, the role of distributed design is more important, especially when communication between different agents is sparse. In addition, by the simulation, we observe that distributed asynchronous algorithm often converges faster than distributed synchronous algorithm. It should be pointed out that although communication time-delay will not make asynchronous algorithm diverge, coordinating the trade-off between communication and computation may further improve the convergence performance of distributed asynchronous algorithm in practice, which is one future research direction of our work.
\begin{table}[!htbp]
	\centering
	\caption{the execution time comparisons}\label{com_table}
	\begin{tabular}{|c|c|c|c|c|c|}
		\hline
		${\rm dimension}$ &${\rm algorithm}$&$sn$&${\rm iterations}$&${\rm time(ms)}$&${\rm error}$\\
		\hline
		\multirow{3}*{8} &$Mixing$&0&291&20&0.00023\\
		\cline{2-6}
		&DSA&5&150&5&0.00023\\
		\cline{2-6}
		&DAA&5&150&4&0.00023\\
		\hline
		\multirow{3}*{18} &$Mixing$ &0&600&41&0.0058\\
		\cline{2-6}
		&DSA&8&390&16&0.0058\\
		\cline{2-6}
		&DAA&8&361&10&0.0058\\
		\hline
	\end{tabular}
\end{table}
\begin{figure}
	\centering
	\includegraphics[width=6cm]{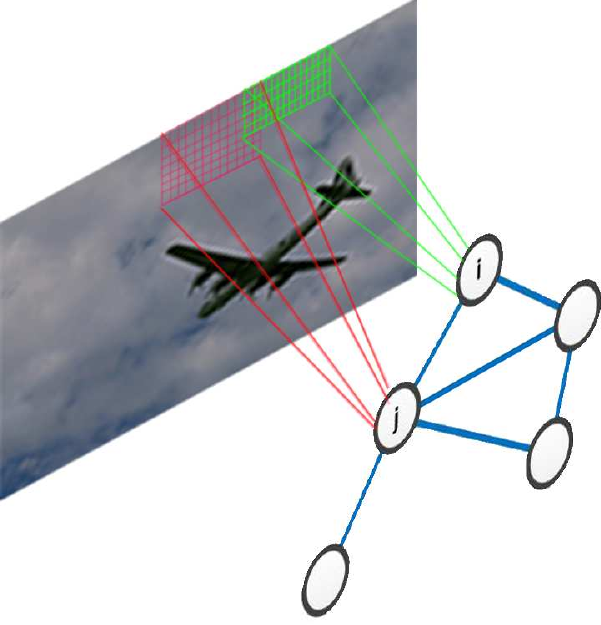}\\
	\caption{Distributed image segmentation set-up. Each agent only has access to a subset (colored grids) of the whole image pixels.}\label{image_cut_show}
	\centering
\end{figure}
\begin{figure}
	\centering
	\includegraphics[width=6cm]{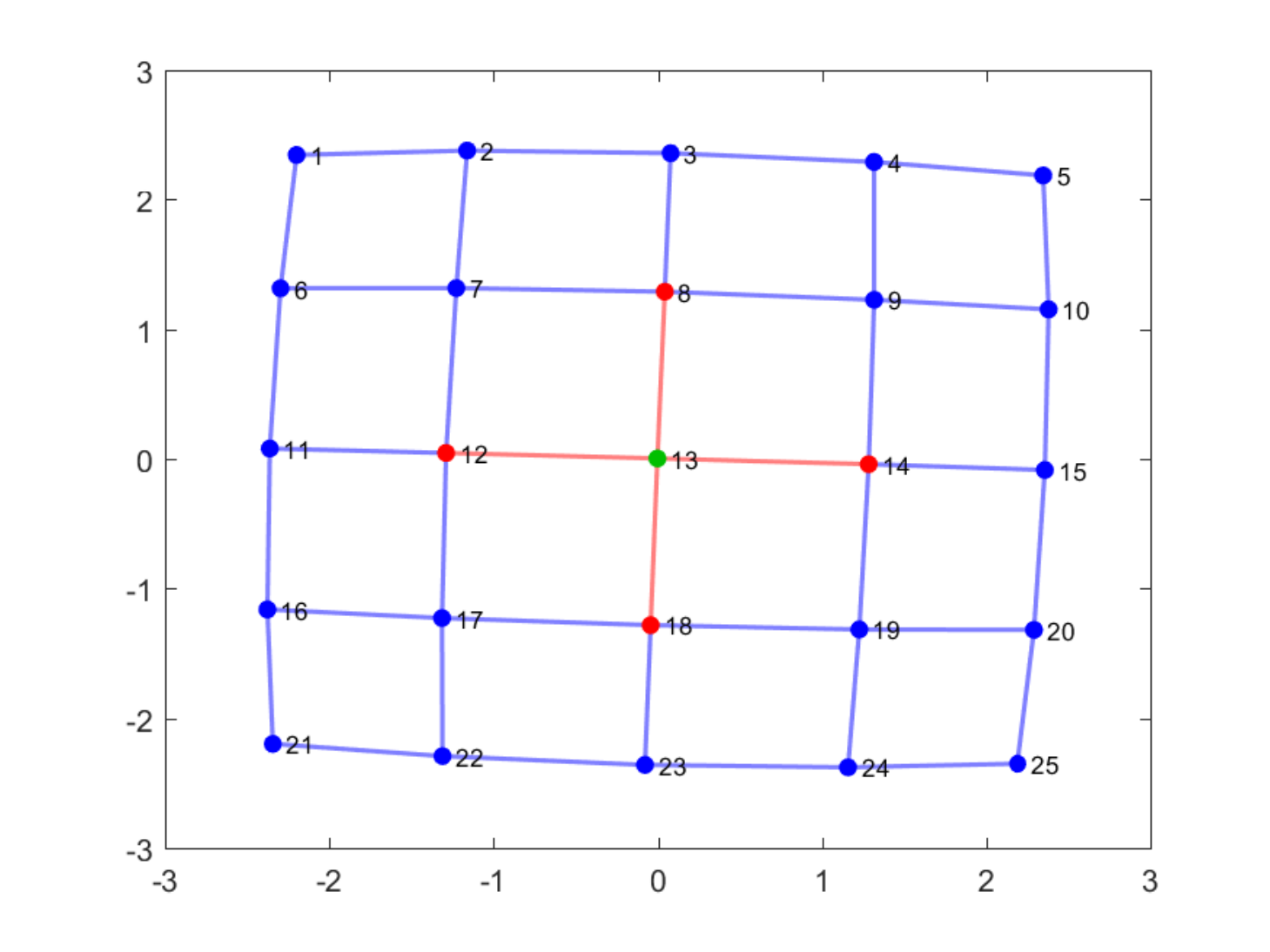}\\
	\caption{The four-connected neighborhood of image segmentation set-up.}\label{pixel_show}
	\centering
\end{figure}
\par \textit{Example 2:} We apply the proposed distributed asynchronous algorithm to solving MAXCUT problems from image segmentation over a multi-agent system, as shown in Fig. \ref{image_cut_show}. There is one edge between agents if there exists an intersection between image pixels of different local images. Hence, there exist several parents for one agent, which is different from the first example, where each agent has only one parent. In this example, we will show that the proposed distributed algorithm efficiently achieves image segmentation.
\par There have been some works applying general graph cut algorithms to image segmentation \cite{cut_survey,cut_2008,image_seg2013,graph_cut_06}. For image segmentation, we need to create a graph representation of the image. One algorithm is to consider assigning each pixel of the image as a node and using a four-connected neighborhood to create the edges\cite{image_seg2013}, as shown in Fig. \ref{pixel_show}. We here only utilize the intensity components of the RGB of all pixels to provide one simple connected matrix $M$, whose $(i,j)$th element related to nodes $(i,j)$ is defined by the following equation \cite{image_seg2013}
\begin{align*}
&M_{i,j}=\\
&{\rm max}((2[\|rgb(i)-rgb(j)\|_2>t]-1)\|rgb(i)-rgb(j)\|_2,0),
\end{align*}
where $rgb(i)$ is the intensity vector of RGB of the $i$th pixel, $t$ is adjustable threshold value. The output of operator ${\rm max}(a,0)$ is the bigger one of $a$ and $0$. Then, the generated matrix $M$ is a typical large-scale sparse matrix. For some algorithms which add seeds to different regions, the only change is the development of matrix elements. We only use the simplest RGB information between different pixels to segment image. However, it should be noted that the proposed algorithms are applicable for general MAXCUT problems (\ref{sdp_pro}) that include more involved development of coefficient matrix elements. In some intelligent algorithms, the graph cut problem is often used as an important pretreatment\cite{contour_p}. Therefore, the large-scale sparse graph cut problem is vital in image segmentation.
\par
\begin{figure}[htbp]
	\centering
	\subfigure[Airplane segmentation]{
		\includegraphics[width=4cm]{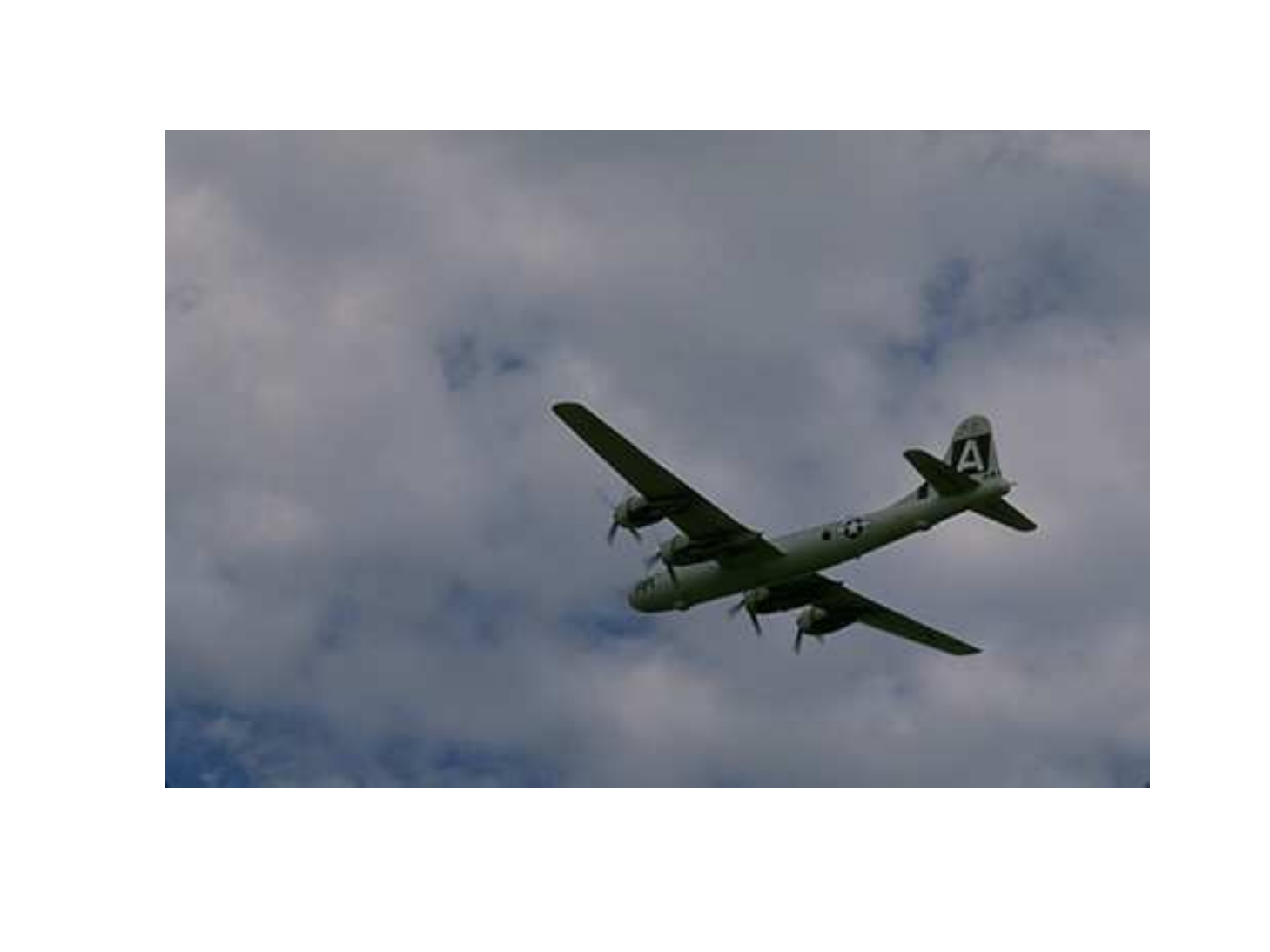}
		\includegraphics[width=4cm]{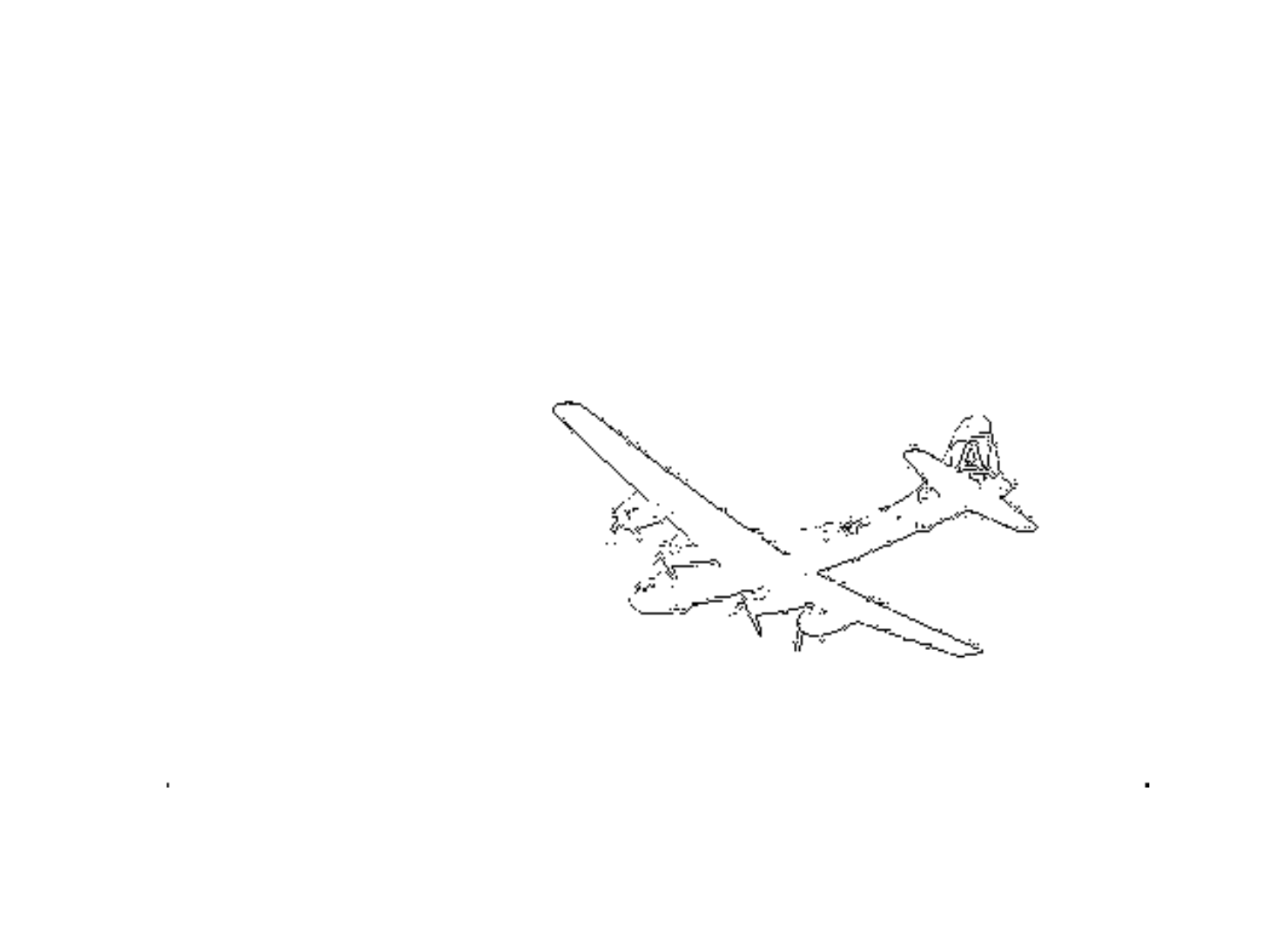}
	}
	
	\subfigure[Church segmentation]{
		
		\includegraphics[width=4cm]{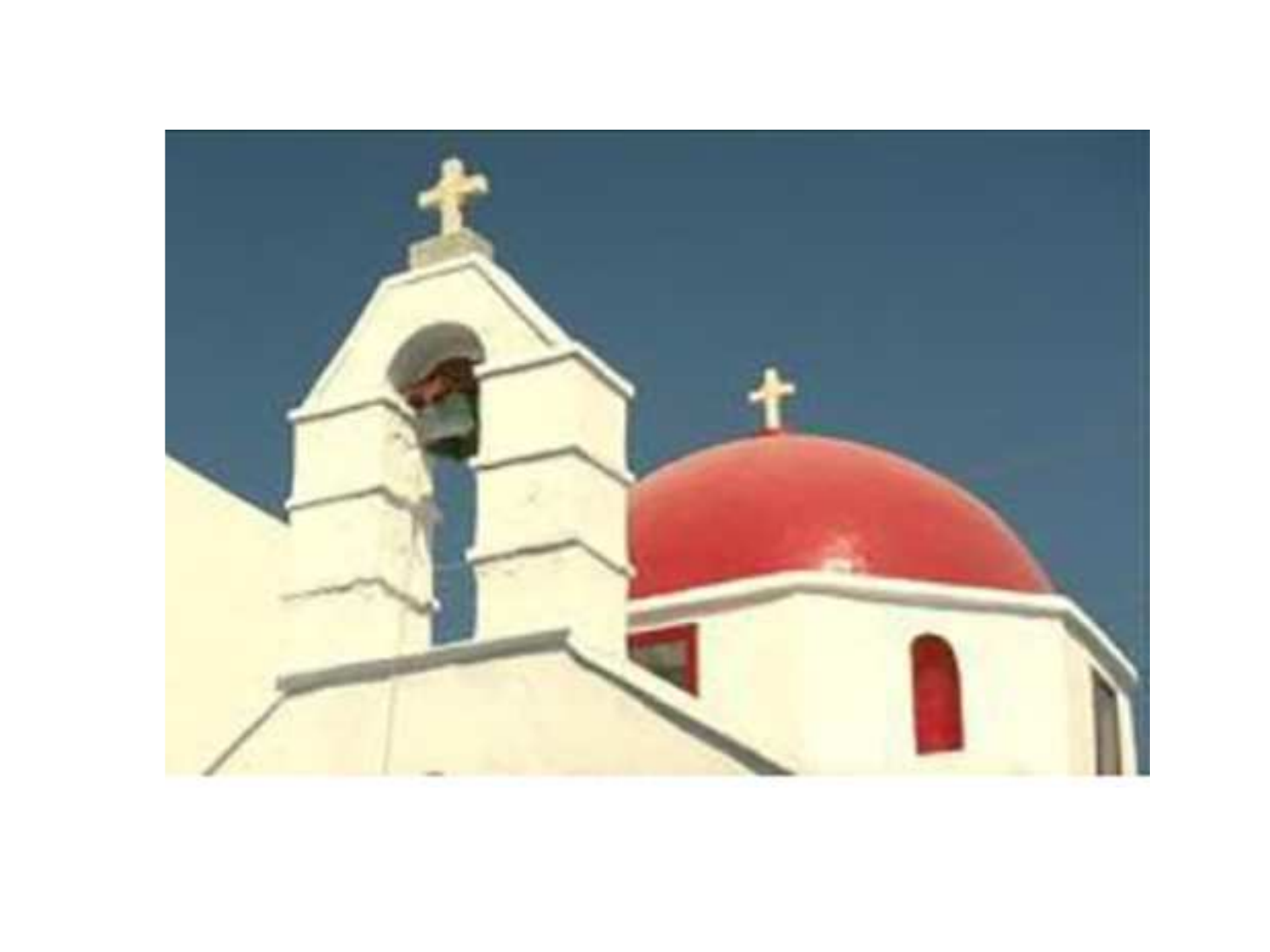}
		\includegraphics[width=4cm]{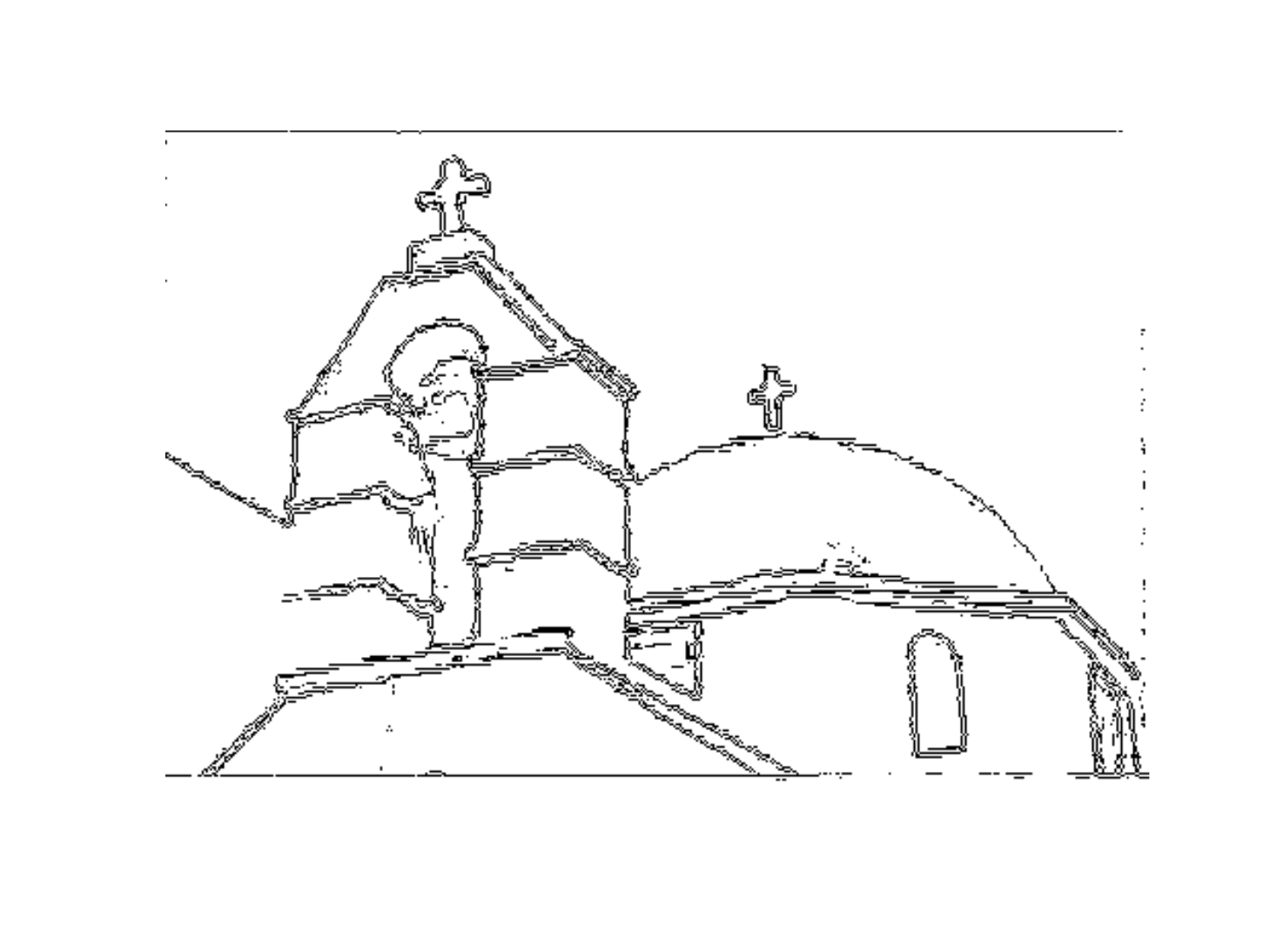}
	}
	\subfigure[Brid segmentation]{
		\includegraphics[width=4cm]{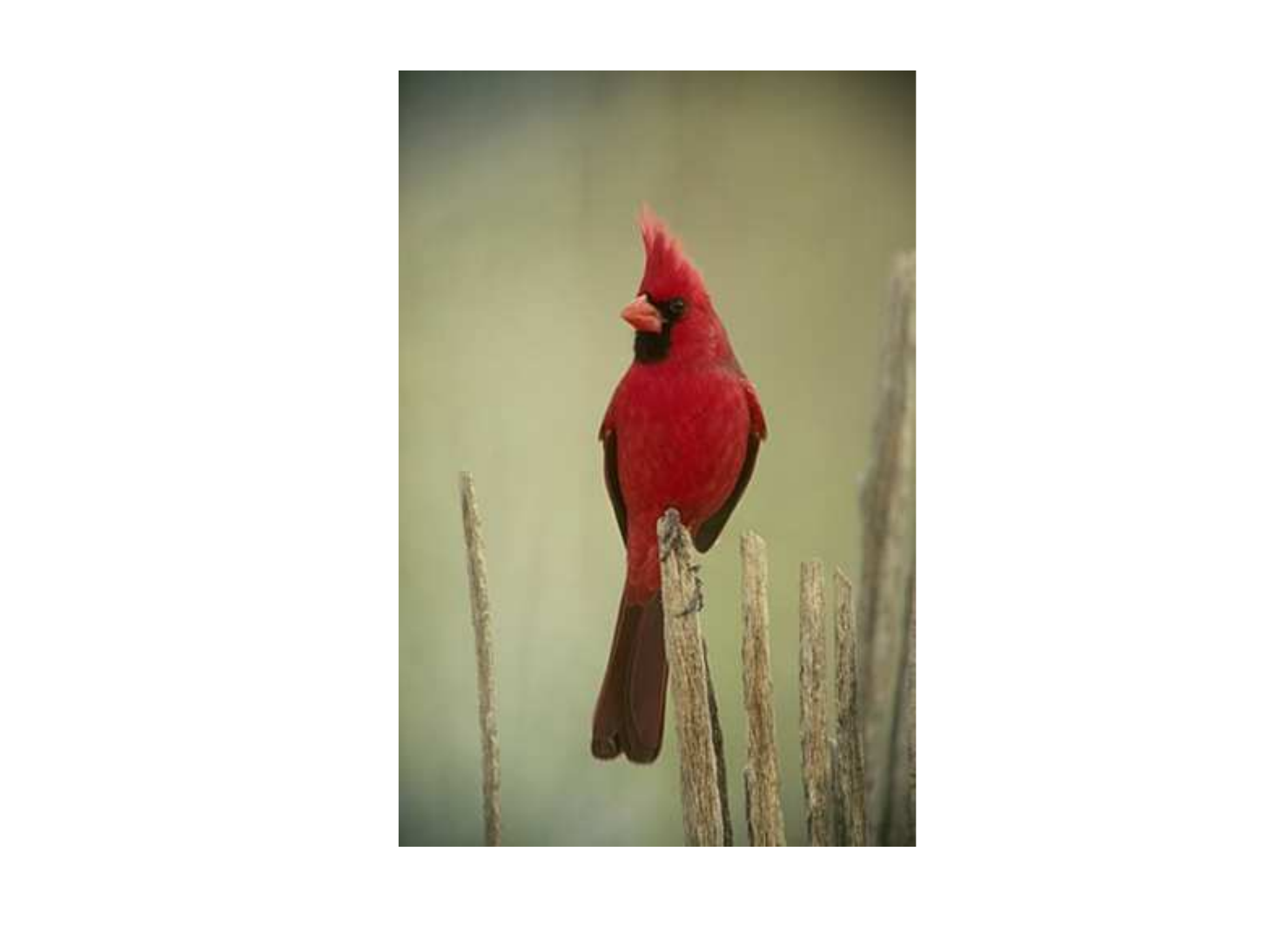}
		\includegraphics[width=4cm]{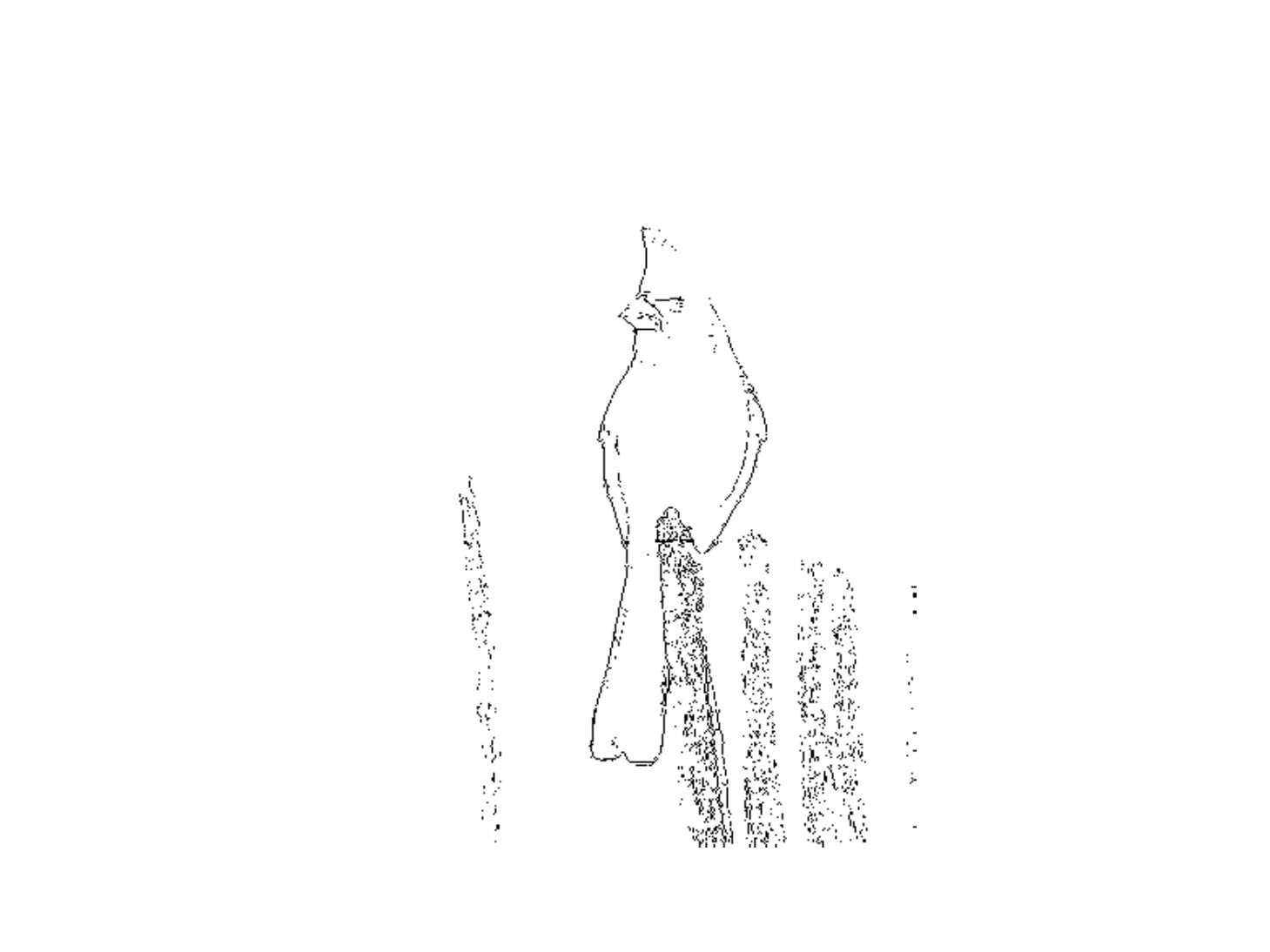}
	}
	
	\caption{Segmentation results of images from Berkeley database}
	\label{image}
\end{figure}

 We apply the proposed distributed asynchronous algorithm on images of the Berkeley database \cite{database}. We have computed the results for three images (Airplane, Church, Bird) in Figure \ref{image}. It is seen that the proposed distributed algorithm achieves image segmentation efficiently. While the existing centralized algorithms can not deal with image segmentation because of the large dimension of image data.

\section{Conclusion}\label{conclusion}
This paper has studied distributed synchronous and asynchronous algorithms for solving large-scale SDP with diagonal constraints by making use of the inherent sparsity of programming and low-rank property of solutions. Each agent updates its local variables by local information and communicating messages over the underlying topology. To handle the communication delays in networks, one distributed asynchronous algorithm is proposed without  global clocks. Although the transformed optimization problem is non-convex, variables of distributed synchronous and asynchronous algorithms eventually converge to optimal solutions and critical points of SDP, repectively.  The efficiency of proposed distributed algorithms is verified by the numerical simulations.
\par Future work involves developing and analyzing communication-efficient distributed algorithms, which balance the computational and communication cost of different agents, for SDP with diagonal constraints. The objective SDP problem of this paper has diagonal constraints. In future, we will further attempt to extend the distributed algorithms to more general semi-definite programs with linear constraints. 
\bibliographystyle{ieeetran}
\bibliography{refer}

\begin{thebibliography}{10}
\providecommand{\url}[1]{#1}
\csname url@samestyle\endcsname
\providecommand{\newblock}{\relax}
\providecommand{\bibinfo}[2]{#2}
\providecommand{\BIBentrySTDinterwordspacing}{\spaceskip=0pt\relax}
\providecommand{\BIBentryALTinterwordstretchfactor}{4}
\providecommand{\BIBentryALTinterwordspacing}{\spaceskip=\fontdimen2\font plus
\BIBentryALTinterwordstretchfactor\fontdimen3\font minus
  \fontdimen4\font\relax}
\providecommand{\BIBforeignlanguage}[2]{{%
\expandafter\ifx\csname l@#1\endcsname\relax
\typeout{** WARNING: IEEEtran.bst: No hyphenation pattern has been}%
\typeout{** loaded for the language `#1'. Using the pattern for}%
\typeout{** the default language instead.}%
\else
\language=\csname l@#1\endcsname
\fi
#2}}
\providecommand{\BIBdecl}{\relax}
\BIBdecl

\bibitem{smart_admm}
E.~{Dall'Anese}, H.~{Zhu}, and G.~B. {Giannakis}, ``Distributed optimal power
  flow for smart microgrids,'' \emph{IEEE Transactions on Smart Grid}, vol.~4,
  no.~3, pp. 1464--1475, 2013.

\bibitem{semi_observer}
F.~R. {Segundo Sevilla}, I.~M. {Jaimoukha}, B.~{Chaudhuri}, and P.~{Korba}, ``A
  semidefinite relaxation procedure for fault-tolerant observer design,''
  \emph{IEEE Transactions on Automatic Control}, vol.~60, no.~12, pp.
  3332--3337, 2015.

\bibitem{2003Semidefinite}
P.~A. Parrilo and S.~Lall, ``Semidefinite programming relaxations and algebraic
  optimization in control,'' \emph{European Journal of Control}, vol.~9, no.~2,
  pp. 307--321, 2003.

\bibitem{center_sdp}
\BIBentryALTinterwordspacing
S.~J. Benson, Y.~Ye, and X.~Zhang, ``Solving large-scale sparse semidefinite
  programs for combinatorial optimization,'' \emph{SIAM Journal on
  Optimization}, vol.~10, no.~2, pp. 443--461, 2000. [Online]. Available:
  \url{https://doi.org/10.1137/S1052623497328008}
\BIBentrySTDinterwordspacing

\bibitem{combina}
Y.~T. {Lee}, A.~{Sidford}, and S.~C. {Wong}, ``A faster cutting plane method
  and its implications for combinatorial and convex optimization,'' in
  \emph{2015 IEEE 56th Annual Symposium on Foundations of Computer Science},
  2015, pp. 1049--1065.

\bibitem{max_like}
M.~{Kim}, J.~{Park}, K.~{Kim}, and J.~{Kim}, ``Exact {ML} criterion based on
  semidefinite relaxation for {MIMO} systems,'' \emph{IEEE Signal Processing
  Letters}, vol.~21, no.~3, pp. 343--346, 2014.

\bibitem{eco_dispatch}
R.~A. {Jabr}, ``Solution to economic dispatching with disjoint feasible regions
  via semidefinite programming,'' \emph{IEEE Transactions on Power Systems},
  vol.~27, no.~1, pp. 572--573, 2012.

\bibitem{matrix_learning}
E.~L. {Hu} and J.~T. {Kwok}, ``Low-rank matrix learning using biconvex
  surrogate minimization,'' \emph{IEEE Transactions on Neural Networks and
  Learning Systems}, vol.~30, no.~11, pp. 3517--3527, 2019.

\bibitem{kernel_matrixlearning}
E.~{Hu}, S.~{Chen}, D.~{Zhang}, and X.~{Yin}, ``Semisupervised kernel matrix
  learning by kernel propagation,'' \emph{IEEE Transactions on Neural
  Networks}, vol.~21, no.~11, pp. 1831--1841, 2010.

\bibitem{distance_matrix}
C.~{Shen}, J.~{Kim}, and L.~{Wang}, ``Scalable large-margin {M}ahalanobis
  distance metric learning,'' \emph{IEEE Transactions on Neural Networks},
  vol.~21, no.~9, pp. 1524--1530, 2010.

\bibitem{1995Improved}
M.~X. Goemans and D.~P. Williamson, ``Improved approximation algorithms for
  maximum cut and satisfiability problems using semidefinite programming,''
  \emph{Journal of the ACM}, vol.~42, no.~6, pp. 1115--1145, 1995.

\bibitem{Waldspurger2012Phase}
I.~Waldspurger, A.~d’Aspremont, and S.~Mallat, ``Phase recovery, maxcut and
  complex semidefinite programming,'' \emph{Mathematical Programming}, vol.
  149, no. 1-2, pp. 47--81, 2012.

\bibitem{2014Exact}
E.~Abbe, A.~S. Bandeira, and G.~Hall, ``Exact recovery in the stochastic block
  model,'' \emph{IEEE Transactions on Information Theory}, vol.~62, no.~1, pp.
  471--487, 2014.

\bibitem{2007Implementation}
B.~Borchers and J.~G. Young, ``Implementation of a primal–dual method for
  {SDP} on a shared memory parallel architecture,'' \emph{Computational
  Optimization and Applications}, vol.~37, no.~3, pp. 355--369, 2007.

\bibitem{BM_smooth}
N.~Boumal, V.~Voroninski, and A.~S. Bandeira, ``The non-convex
  {B}urer–{M}onteiro approach works on smooth semidefinite programs,'' in
  \emph{Proceedings of the 30th International Conference on Neural Information
  Processing Systems}, ser. NIPS’16.\hskip 1em plus 0.5em minus 0.4em\relax
  Red Hook, NY, USA: Curran Associates Inc., 2016, p. 2765–2773.

\bibitem{2015Phase}
E.~J. Candes, X.~Li, and M.~Soltanolkotabi, ``Phase retrieval via {W}irtinger
  flow: Theory and algorithms,'' \emph{IEEE Transactions on Information
  Theory}, vol.~61, no.~4, 2015.

\bibitem{2017Solving}
S.~Mei, T.~Misiakiewicz, A.~Montanari, and R.~I. Oliveira, ``Solving {SDP}s for
  synchronization and maxcut problems via the grothendieck inequality,''
  \emph{arXiv:Optimization and Control}, 2017.

\bibitem{wang2017mixing}
P.-W. Wang, W.-C. Chang, and J.~Z. Kolter, ``The {M}ixing method: low-rank
  coordinate descent for semidefinite programming with diagonal constraints,''
  \emph{arXiv:Optimization and Control}, 2017.

\bibitem{Zhang2012}
M.~{Zhang}, ``A second order mehrotra-type predictor-corrector algorithm for
  semidefinite optimization,'' \emph{Journal of Systems Science and
  Complexity}, vol.~25, no.~6, pp. 1108--1121, 2012.

\bibitem{BM_SPLR}
S.~{Burer} and R.~D. {Monteiro}, ``A nonlinear programming algorithm for
  solving semidefinite programs via low-rank factorization,''
  \emph{Mathematical Programming}, vol.~95, no.~2, pp. 329--357, 2 2003.

\bibitem{lowrank_sdp}
\BIBentryALTinterwordspacing
M.~Journée, F.~Bach, P.-A. Absil, and R.~Sepulchre, ``Low-rank optimization on
  the cone of positive semidefinite matrices,'' \emph{SIAM Journal on
  Optimization}, vol.~20, no.~5, pp. 2327--2351, 2010. [Online]. Available:
  \url{https://doi.org/10.1137/080731359}
\BIBentrySTDinterwordspacing

\bibitem{erdogdu2018convergence}
M.~A. Erdogdu, A.~Ozdaglar, P.~A. Parrilo, and N.~D. Vanli, ``Convergence rate
  of block-coordinate maximization {B}urer-{M}onteiro method for solving large
  {SDP}s,'' \emph{arXiv:Optimization and Control}, 2018.

\bibitem{Sylvester_matrix}
W.~{Deng}, X.~{Zeng}, and Y.~{Hong}, ``Distributed computation for solving the
  sylvester equation based on optimization,'' \emph{IEEE Control Systems
  Letters}, vol.~4, no.~2, pp. 414--419, 2020.

\bibitem{linear_eq}
G.~{Shi}, B.~D.~O. {Anderson}, and U.~{Helmke}, ``Network flows that solve
  linear equations,'' \emph{IEEE Transactions on Automatic Control}, vol.~62,
  no.~6, pp. 2659--2674, 2017.

\bibitem{Deng2019NetworkFT}
W.~Deng, Y.~Hong, B.~Anderson, and G.~Shi, ``Network flows that solve sylvester
  matrix equations,'' \emph{arXiv: Optimization and Control}, 2019.

\bibitem{2016Implementing}
J.~Yang, X.~Meng, and M.~W. Mahoney, ``Implementing randomized matrix
  algorithms in parallel and distributed environments,'' \emph{Proceedings of
  the IEEE}, vol. 104, no.~1, pp. 58--92, 2016.

\bibitem{matrix_equation}
W.~Li, X.~Zeng, Y.~Hong, and J.~Haibo, ``Distributed design for nuclear norm
  minimization of linear matrix equation with constraints,'' \emph{IEEE
  Transactions on Automatic Control}, pp. 1--1, 2020.

\bibitem{Lyapunov_matrix}
X.~Jiang, X.~Zeng, J.~Sun, and J.~Chen, ``Distributed solver for discrete-time
  {L}yapunov equations over dynamic networks with linear convergence rate,''
  \emph{IEEE Transactions on Cybernetics}, pp. 1--10, 2020.

\bibitem{Hong2016}
Z.~{Deng} and Y.~{Hong}, ``Multi-agent optimization design for autonomous
  lagrangian systems,'' \emph{Unmanned Systems}, vol.~4, no.~1, pp. 5--13,
  2016.

\bibitem{consensus_shi}
A.~{Fontan}, G.~{Shi}, X.~{Hu}, and C.~{Altafini}, ``Interval consensus for
  multiagent networks,'' \emph{IEEE Transactions on Automatic Control},
  vol.~65, no.~5, pp. 1855--1869, 2020.

\bibitem{liang_opti}
S.~{Liang}, L.~Y. {Wang}, and G.~{Yin}, ``Distributed smooth convex
  optimization with coupled constraints,'' \emph{IEEE Transactions on Automatic
  Control}, vol.~65, no.~1, pp. 347--353, 2020.

\bibitem{Wang2020}
K.~{Wang}, Z.~{Fu}, Q.~{Xu}, D.~{Chen}, L.~{Wang}, and W.~{Yu}, ``Distributed
  fixed step-size algorithm for dynamic economic dispatch with power flow
  limits,'' \emph{Science China Information Sciences}, vol.~64, no.~1, p.
  112202, 2020.

\bibitem{ADMM_CDC}
R.~{Madani}, A.~{Kalbat}, and J.~{Lavaei}, ``{ADMM} for sparse semidefinite
  programming with applications to optimal power flow problem,'' in \emph{2015
  54th IEEE Conference on Decision and Control (CDC)}, Osaka, 2015, pp.
  5932--5939.

\bibitem{fast_dis}
A.~{Kalbat} and J.~{Lavaei}, ``A fast distributed algorithm for decomposable
  semidefinite programs,'' in \emph{2015 54th IEEE Conference on Decision and
  Control (CDC)}, Osaka, 2015, pp. 1742--1749.

\bibitem{semi_power}
H.~{Zhu} and G.~B. {Giannakis}, ``Power system nonlinear state estimation using
  distributed semidefinite programming,'' \emph{IEEE Journal of Selected Topics
  in Signal Processing}, vol.~8, no.~6, pp. 1039--1050, 2014.

\bibitem{asyn_admm}
C.~{Chang}, J.~{Cortés}, and S.~{Martínez}, ``Scheduled-asynchronous
  distributed algorithm for optimal power flow,'' \emph{IEEE Transactions on
  Control of Network Systems}, vol.~6, no.~1, pp. 261--275, 2019.

\bibitem{dis_pd}
S.~K. Pakazad, A.~Hansson, M.~S. Andersen, and I.~Nielsen, ``Distributed
  primal–dual interior-point methods for solving tree-structured coupled
  convex problems using message-passing,'' \emph{Optimization Methods and
  Software}, vol.~32, no.~3, pp. 401--435, 2017.

\bibitem{ricatti_zeng}
X.~{Zeng}, J.~{Chen}, and Y.~{Hong}, ``Distributed optimization design for
  computation of algebraic {R}iccati inequalities,'' \emph{IEEE Transactions on
  Cybernetics}, pp. 1--12, 2020.

\bibitem{2016community}
A.~S. Bandeira, N.~Boumal, and V.~Voroninski, ``On the low-rank approach for
  semidefinite programs arising in synchronization and community detection,''
  in \emph{29th Annual Conference on Learning Theory}, ser. Proceedings of
  Machine Learning Research, V.~Feldman, A.~Rakhlin, and O.~Shamir, Eds.,
  vol.~49.\hskip 1em plus 0.5em minus 0.4em\relax Columbia University, New
  York, New York, USA: PMLR, 23--26 Jun 2016, pp. 361--382.

\bibitem{mincut_graph}
Y.~{Boykov} and V.~{Kolmogorov}, ``An experimental comparison of
  min-cut/max-flow algorithms for energy minimization in vision,'' \emph{IEEE
  Transactions on Pattern Analysis and Machine Intelligence}, vol.~26, no.~9,
  pp. 1124--1137, 2004.

\bibitem{chor_spar}
R.~P. {Mason} and A.~{Papachristodoulou}, ``Chordal sparsity, decomposing
  {SDP}s and the {L}yapunov equation,'' in \emph{2014 American Control
  Conference}, Portland, USA, 2014, pp. 531--537.

\bibitem{SDP_rank}
G.~Pataki, ``On the rank of extreme matrices in semidefinite programs and the
  multiplicity of optimal eigenvalues,'' \emph{Mathematics of Operations
  Research}, vol.~23, no.~2, pp. 339--358, 1998.

\bibitem{saddle_escape}
J.~D. {Lee}, I.~{Panageas}, G.~{Piliouras}, M.~{Simchowitz}, M.~I. {Jordan},
  and B.~{Recht}, ``First-order methods almost always avoid strict saddle
  points,'' \emph{Mathematical Programming}, vol. 176, no.~1, pp. 311--337, 7
  2019.

\bibitem{inci2012regularity}
H.~Inci, T.~Kappeler, and P.~Topalov, ``On the regularity of the composition of
  diffeomorphisms,'' \emph{arXiv:Analysis of PDEs}, 2012.

\bibitem{shub_book}
M.~Shub, \emph{Global stability of dynamical systems}.\hskip 1em plus 0.5em
  minus 0.4em\relax Springer Science $\&$ Business Media, 2013.

\bibitem{paral_distri_book}
D.~P. Bertsekas and J.~N. Tsitsiklis, \emph{Parallel and distributed
  computation: numerical methods}.\hskip 1em plus 0.5em minus 0.4em\relax
  Belmont Massachusetts: Athena Scientific, 1997.

\bibitem{distri_sed_tac}
S.~K. {Pakazad}, A.~{Hansson}, M.~S. {Andersen}, and A.~{Rantzer},
  ``Distributed semidefinite programming with application to large-scale system
  analysis,'' \emph{IEEE Transactions on Automatic Control}, vol.~63, no.~4,
  pp. 1045--1058, 2018.

\bibitem{cut_survey}
F.~{Yi} and I.~{Moon}, ``Image segmentation: A survey of graph-cut methods,''
  in \emph{2012 International Conference on Systems and Informatics
  (ICSAI2012)}, Yantai, 2012, pp. 1936--1941.

\bibitem{cut_2008}
S.~{Vicente}, V.~{Kolmogorov}, and C.~{Rother}, ``Graph cut based image
  segmentation with connectivity priors,'' in \emph{2008 IEEE Conference on
  Computer Vision and Pattern Recognition}, Anchorage, AK, 2008, pp. 1--8.

\bibitem{image_seg2013}
S.~{de Sousa}, Y.~{Haxhimusa}, and W.~{Kropatsch}, ``Estimation of distribution
  algorithm for the {M}ax-{C}ut problem,'' in \emph{Graph-Based Representations
  in Pattern Recognition (GbRPR 2013)}, vol. 7877, Berlin, Heidelberg, 2013,
  pp. 244--253.

\bibitem{graph_cut_06}
Y.~{Boykov} and G.~{Funka-Lea}, ``Graph cuts and efficient {ND} image
  segmentation,'' \emph{International Journal of Computer Vision}, vol.~70,
  no.~2, pp. 109--131, 2006.

\bibitem{contour_p}
P.~{Arbelaez}, M.~{Maire}, C.~{Fowlkes}, and J.~{Malik}, ``Contour detection
  and hierarchical image segmentation,'' \emph{IEEE Transactions on Pattern
  Analysis and Machine Intelligence}, vol.~33, no.~5, pp. 898--916, 2011.

\bibitem{database}
D.~{Martin}, C.~{Fowlkes}, D.~{Tal}, and J.~{Malik}, ``A database of human
  segmented natural images and its application to evaluating segmentation
  algorithms and measuring ecological statistics,'' in \emph{Proceedings Eighth
  IEEE International Conference on Computer Vision. ICCV 2001}, vol.~2, British
  Columbia, Canada, 2001, pp. 416--423.

\end{thebibliography}

\end{document}